\documentclass[twoside,11pt]{article}

\usepackage{jmlr2e}

\usepackage{amsmath,amssymb, makecell, listings, mathtools, lsc, booktabs,colortbl,bibentry,paralist}
\usepackage{ mathrsfs,graphics,epsfig}
\usepackage[font=small, labelfont=bf]{caption}
\usepackage[font=footnotesize, labelfont=bf]{subcaption}
\usepackage{setspace}
\usepackage{pst-node}
\usepackage{tikz-cd} 

\usepackage{todonotes}
\newcounter{todocounter}

\presetkeys{todonotes}{inline, color=yellow!30}{}
\usepackage{nameref}

\usepackage{algorithm}[chapter]
\usepackage{algorithmic}
\floatname{algorithm}{Algorithm}

\newcommand{\R}{\ensuremath{\mathbb{R}}}

\newcommand{\Z}{\ensuremath{\mathbb{Z}}}

\DeclareMathOperator{\Ima}{im}

\jmlrheading{1}{}{1-48}{4/00}{10/00}{meila00a}{Bernadette J. Stolz}

\ShortHeadings{Outlier-robust subsampling techniques for PH}{Stolz}
\firstpageno{1}

\begin{document}

\title{Outlier-robust subsampling techniques for persistent homology}

\author{\name Bernadette J. Stolz \email stolz@maths.ox.ac.uk \\
       \addr Mathematical Institute\\
        University of Oxford, Woodstock Rd\\
       Oxford OX2 6GG, United Kingdom}

\editor{} 

\maketitle

\begin{abstract}
In recent years, persistent homology (PH) has been successfully applied to real-world data in many different settings.
Despite significant computational advances, PH algorithms do not yet scale to large datasets preventing interesting applications.
One approach to address computational issues posed by PH is to select a set of landmarks by subsampling from the data. 
Currently, these landmark points are chosen either at random or using the maxmin algorithm. 
Neither is ideal as random selection tends to favour dense areas of the data while the  maxmin algorithm is very sensitive to noise. Here, we propose a novel approach to select landmarks specifically for PH that preserves coarse topological information of the original dataset. 
Our method is motivated by the Mayer-Vietoris sequence and requires only local PH computation thus enabling efficient computation. We test our landmarks on artificial datasets which contain different levels of noise and compare them to standard landmark selection techniques. 
We demonstrate that our landmark selection outperforms standard methods as well as a subsampling technique based on an outlier-robust version of the $k$--means algorithm for low sampling densities in noisy data with respect to robustness to outliers.
\end{abstract}


\section{Introduction}

One of the practical challenges of applying persistent homology (PH) is that it is computationally difficult for large datasets. 
Building a filtration on a dataset with $N$ points results in spaces of the size $\mathcal{O}(2^N)$, although in practice this can be reduced to $\mathcal{O}(N^{\hat{n}+1})$ by posing a limit $\hat{n}$ on the dimension of the topological features considered~\citep{Otter2017}.
Following the construction of a filtration, the algorithm for computing PH in the worst case has a complexity of $\mathcal{O}(\mathtt{k}^3)$, where $\mathtt{k}$ is the number of points, edges, triangles, and higher-dimensional connections constructed from the data points by the filtration (the number of such connections can be up to $2^N-1$), although the complexity is often linear in practice~\citep{Otter2017,Edelsbrunner2002}. Methods exist to approximate specific filtrations or to reduce the sizes of the vector spaces associated to the data by a filtration, for an overview of such methods that have been
 implemented in software packages, see~\citep{Otter2017}. Despite such improvements, computation is still very challenging on large datasets and can pose a hard limit on the filtrations that can be applied to a particular dataset. 
In such cases it can become necessary to preprocess data before applying PH.
For point cloud data one can, for example, use subsampling techniques to identify so-called \emph{landmarks} of the dataset and then define a filtration on the landmarks~\citep{deSilva2004}. A filtration on well chosen landmarks retains topologically important global information about the full dataset and one can even choose to include information from non-landmark points when constructing the filtration.
An example for such a filtration is the \emph{lazy witness filtration} which was first introduced by \citep{deSilva2004} and has been used to study noisy artificial datasets \citep{Kovacev2012}, primary visual cortex cell populations ~\citep{Singh2008}, and cancer gene expression data ~\citep{Lockwood2015}. Roughly, the lazy witness filtration consists of the following steps\footnote{For the full definition see~\citep{deSilva2004}.}:
\begin{enumerate}
\item Selection of a (small) subset\footnote{Although there is no systematic lower bound for the number of landmark points, \citep{deSilva2004} suggest using $ > 5\%$ of the data points as landmarks.} of \textit{landmark points} $L$ from the dataset $D$.
\item Construction of a \textit{lazy witness filtration} on the landmarks $L$ where the landmarks are vertices and data points from the full dataset $D$ can serve as \textit{witnesses} for higher order interactions between sets of landmarks in the filtration.
\end{enumerate}
Even though the choice of landmarks from the data inevitably has a large influence on the results that can be obtained, currently there are only two standard approaches to select landmarks: uniform random selection and selection via the \emph{maxmin} algorithm. Neither is ideal and in particular the  maxmin algorithm tends to include outliers \citep{deSilva2004, Adams2014}, which poses problems since large real-world datasets often include noise. In the biological application of the lazy witness filtration in~\citep{Lockwood2015}, for example, the  maxmin algorithm leads to the discovery of loops in the dataset which we could not reproduce using uniform random landmark selection or when discarding a small proportion of the initially chosen  maxmin landmarks from the dataset and choosing a new set of landmarks with the  maxmin algorithm. 
We also note that the results of the lazy witness filtration are difficult to interpret. In addition to the mentioned disadvantages of the maxmin algorithm, neither of the proposed landmark selection methods were designed specifically for PH. While there are other methods that address subsampling for PH~\citep{Cohen-Steiner2007,Niyogi2008,Dufresne2019} these do not explicitly consider noisy data. 
Because existing approaches are not ideal, it is desirable to develop new methods that lead to a reduction of large and noisy data while preserving topological properties. Ideally, the reduced dataset can be used as input for PH directly without additional preprocessing steps.

While one of the appeals of PH is its ability to study multi-scale datasets globally, local PH around a data point can also produce useful insights, see for example ~\citep{Bendich2012,Ahmed2014,Fasy2016,Stolz2020,Wheeler2021}.
Here, we present a novel landmark selection technique that allows us to obtain outlier-robust PH landmarks from noisy point cloud data. Our method is motivated by the Mayer--Vietoris sequence and relies on computing the Vietoris--Rips filtration~\citep{Carlsson2009, Ghrist2008} locally, in a small neighbourhood around each data point, where the aforementioned computational problems of PH vanish. Our algorithm then uses PH output to define a score for each data point which allows us to identify suitable candidates for landmarks. We investigate two different flavours of our approach, one using high scores to identify candidate landmarks and one using low scores.
We apply our approach to very simple artificial datasets that consist of signal points sampled from an object with a topologically interesting structure, such as a sphere, a torus, or a Klein bottle, as well as noise points. 
A consequence of the inclusion of noise is that a large dataset cannot be reduced by applying subsampling techniques developed to infer the (persistent) homology of data~\citep{Cohen-Steiner2007,Niyogi2008,Dufresne2019}.
The (global) PH of landmarks selected in our proposed method is close to the PH of signal points in the original dataset. 
In comparison to existing landmark selection procedures, we demonstrate that our landmarks based on local PH perform very well on our datasets, in particular for low sampling densities, with respect to robustness to outliers. We further investigate the performance of our landmarks in comparison with an outlier-robust version of the $k$--means algorithm~\citep{Chawla2013} which we modify for landmark selection and again find that they are superior for low sampling densities. Our PH landmarks present a valuable addition to the two existing techniques for landmarks selection for PH and we believe that they may also be of interest in a broader data science context.

Local homology has previously been used to infer the stratification of data~\citep{Bendich2012,Mileyko2021,Nanda2020,Robinson2018,Stolz2020} as well as road network analysis \citep{Ahmed2014,Fasy2016} and, more recently, the analysis of topological complexity of data as it passes through different layers of neural networks \citep{Wheeler2021}.
A similar Mayer-Vietoris sequence as the one that motivates our landmark selection approach has been used to show Wasserstein stability bounds of the Vietoris-Rips filtration on finite point clouds \citep{Skraba2020} 
and local computations motivated by the Mayer-Vietoris sequence have also been applied with the goal of parallelising PH computations ~\citep{casas2019}. 
Although the idea behind our notion of local persistent homology is similar to some of those used previously~\citep{Ahmed2014, Bendich2012, Fasy2016,Wheeler2021}, in comparison to these approaches, our definition does not use relative homology and can be readily computed by considering the Vietoris--Rips complex on points in a small local neighbourhood.

Our paper is structured as follows: in Section~\ref{subsec:LandmarksSoFar}, we review existing standard techniques for landmark selection. We then mathematically motivate and introduce our novel landmark selection technique, which we refer to as PH landmarks, and also present an outlier-robust version of the $k$--means algorithm~\citep{Chawla2013} modified for landmark selection, $k--$ landmarks, in Section \ref{subsec:LandmarksNovel}. In Section~\ref{subsec:LandmarkData} we introduce our datasets. We compare our PH landmarks to existing techniques as well as $k--$ landmarks in Section~\ref{subsec:LandmarkResults} followed by a discussion in Section~\ref{subsec:LandmarkDiscussion}.

\section{Existing landmark selection methods}\label{subsec:LandmarksSoFar}

The two standard methods for selecting landmarks $L\subset D$ in a dataset ${D = \{y_1,\dots,y_N\}}$ are random landmark selection and the  maxmin algorithm. Both are implemented as standard procedures for use in combination with the lazy witness filtration in the PH software package {\sc javaPlex}~\citep{Javaplex}.

\subsection{Random landmark selection}

\begin{figure}
\includegraphics[width=.48\textwidth]{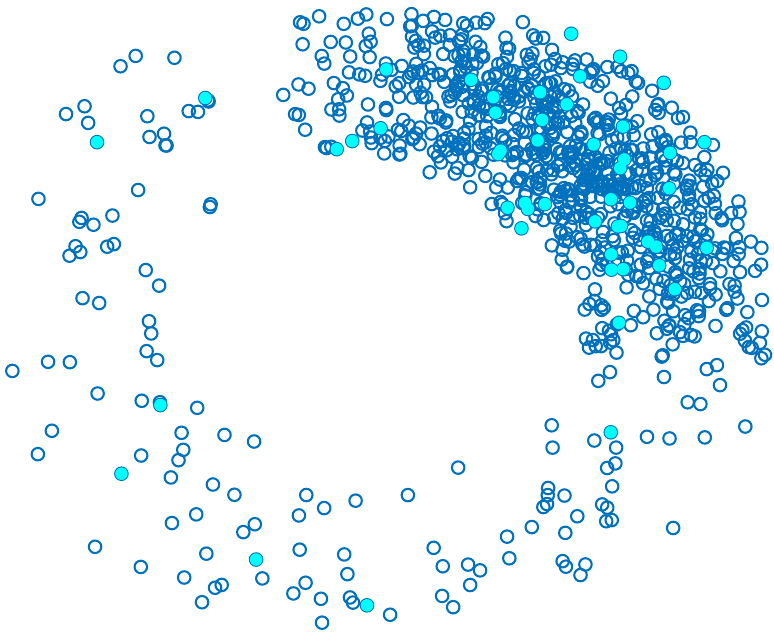}
\includegraphics[width=.48\textwidth]{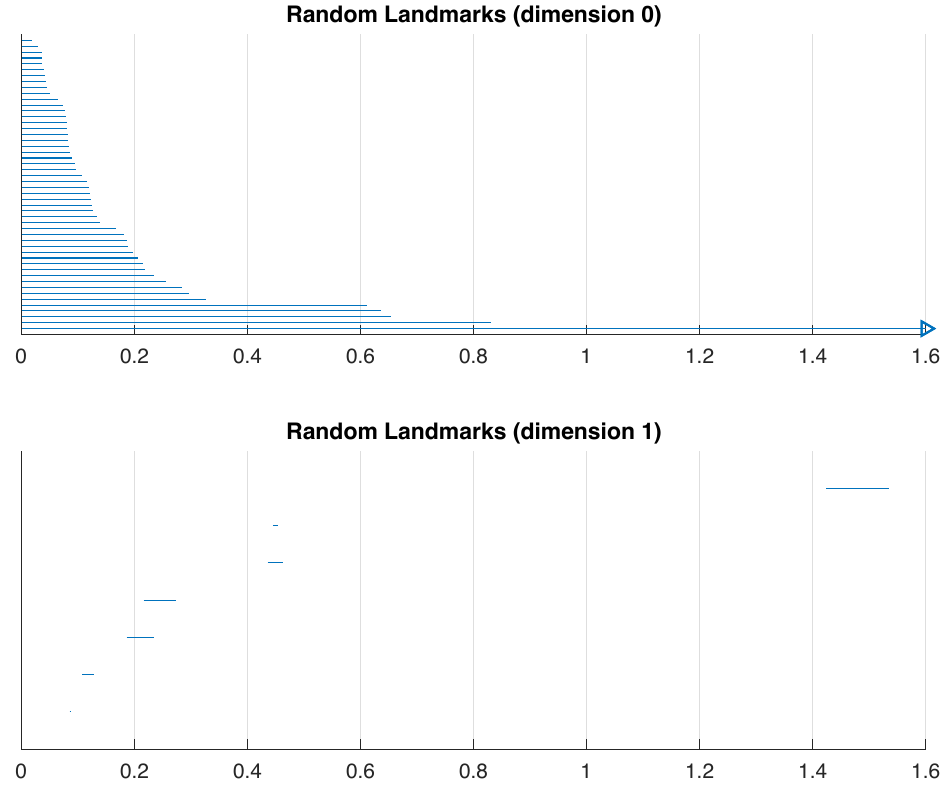}
\caption[Example of a point cloud and landmarks selected at random.]{\label{Fig:RandomFail} Example of a point cloud and landmarks selected at random (landmarks are shown in cyan). We observe that the dimension 1 barcode based on a Vietoris--Rips filtration on the selected landmarks does not capture the persistent homology of the point cloud correctly.}
\end{figure}

The simplest way to choose landmarks $L = \{ l_1, l_2,\dots,l_m\}$ from a point cloud $D$ is to select $m$ points from $D$ uniformly at random. For datasets whose points are evenly distributed, random selection achieves good coverage at a small computational cost. However, as soon as there are large differences in the density of the data, random selection will favour points from more dense regions, which can result in landmarks that do not represent the point cloud well. In extreme cases, the landmarks do not carry any topological similarity to the original point cloud. We show such an example in Fig.~\ref{Fig:RandomFail}.

\subsection{The  maxmin algorithm}

\begin{figure}
\includegraphics[width=.48\textwidth]{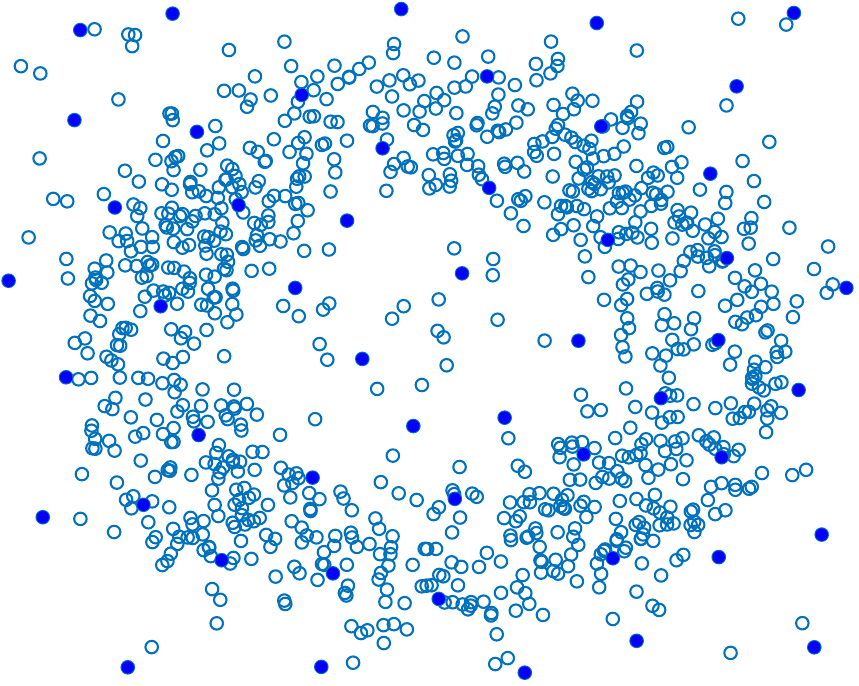}
\includegraphics[width=.48\textwidth]{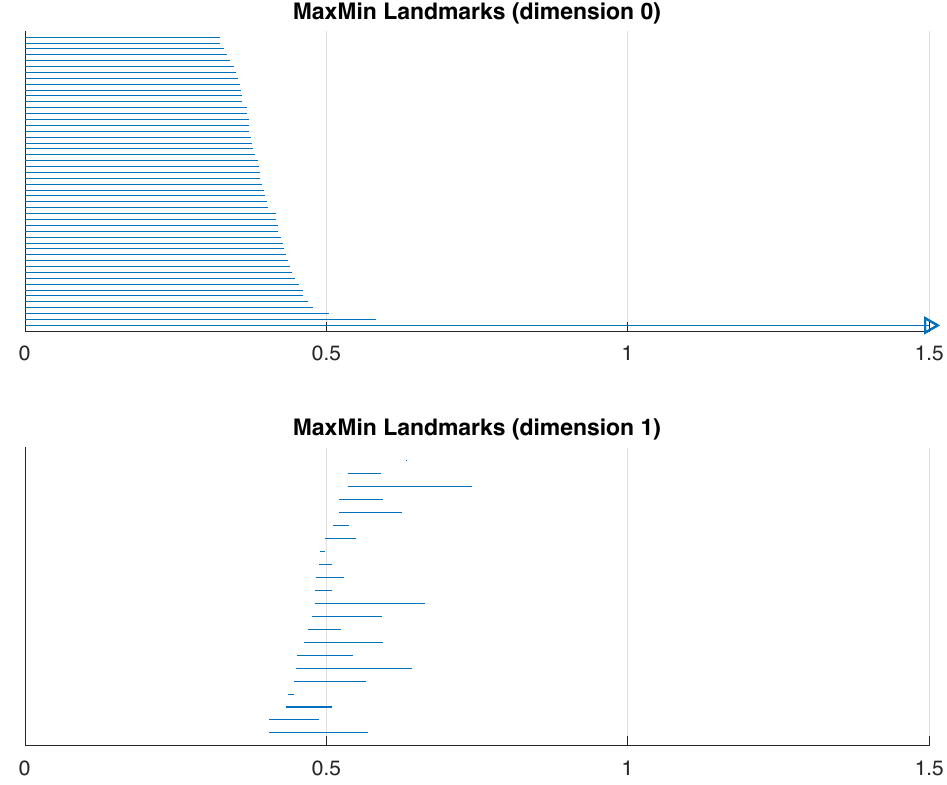}
\caption[Example of a point cloud and landmarks selected by the  maxmin algorithm.]{\label{Fig:MaxMinFail}Example of a point cloud and landmarks selected by the  maxmin algorithm (landmarks are shown in dark blue). We observe that the dimension 1 barcode based on a Vietoris--Rips filtration on the selected landmarks does not capture the PH of the point cloud correctly.}
\end{figure}

The sequential maxmin algorithm chooses the first landmark $l_1 \in D$ randomly. Inductively, for $i \geq 2$ and a landmark set $L_{i-1} = \{ l_1, l_2,\dots,l_{i-1}\}$, the algorithm selects the next landmark $l_i \in D\setminus{L_{i-1} }$ such that for a chosen metric ${d:D\times D \rightarrow \R}$ with  $d(y,L) = \min_{l \in L}{d(y,l)}$ the function mapping $$y \mapsto d(y,L_{i-1}),$$  is maximised for $y \in D$. 
 We show pseudocode for the procedure in Algorithm~\ref{alg:maxmin} in Appendix \ref{sec:Pseudocodemaxmin}.
The method has been used successfully for image data in~\citep{Adams2009,Carlsson2008}.
Landmarks chosen in this way cover the dataset well, are evenly spaced, and represent the underlying topological features better than landmarks selected at random. 
The algorithm does, however, tend to include outliers~\citep{deSilva2004, Adams2014}. We show an example of a point cloud where the selected landmarks do not represent the underlying topology of the data correctly in Fig.~\ref{Fig:MaxMinFail}.

According to \citep{deSilva2004}, the maxmin algorithm is best suited to produce the qualities desired from landmarks. They do not recommend the use of clustering algorithms as an alternative due to the high computational cost and potential to accentuate accidental features. 

\subsection{Dense core subsets}

The tutorial accompanying the PH software package {\sc javaPlex}~\citep{Adams2014} mentions that using so-called dense core subsets before applying the  maxmin algorithm can help overcome the selection of outliers as landmarks. This approach, however, is not considered as one of the standard approaches for landmark selection and we did not find examples where it was used in practice for this purpose. De Silva and Carlsson, 2004 
and Carlsson \emph{et al.}, 2008 
apply dense core subsets to identify dense regions in their datasets. The authors subsequently use the  maxmin landmarks to study the topology of these dense regions.

Dense core subsets are based on assigning density values to every point: for an integer $K$, the density value assigned to a point $y \in D$ is $\frac{1}{\rho_K(y)}$, where $\rho_K(y)$ is the distance to the $K$-th nearest neighbour of $y$. Large values of $K$ provide a measure of the global density around the point in the dataset, while smaller values of $K$ give a more local perspective. Using the density values, one can select the $m$ densest points in the dataset as a dense core subset. Given a dataset, it is not clear what values of $K$ to use. Different values for $K$ and $m$ can produce markedly different subsets~\citep{deSilva2004}.

For our comparisons in Subsection~\ref{Subsubsec:DenseCoreComparison}, instead of selecting a dense subset and then performing the  maxmin algorithm to select landmarks, as proposed by \citep{Adams2014}, we choose the $m$ densest points in the data as landmarks. This enables us to determine how the information on which our own landmark selection technique (PH landmarks) is based differs from that given by the $K$-th nearest neighbour of a data point. 

\section{Proposed landmark selection methods}\label{subsec:LandmarksNovel}

Currently existing methods for landmark selection that enable subsequent PH analysis on large and noisy point cloud data are not ideal and have, in particular, not been designed with PH in mind. 
De Silva and Carlsson, 2004 
state the most pertinent qualities for a landmark set to be good coverage of the dataset and even spacing spacing of the landmarks.
While these are certainly important properties for many datasets,
we find that in the context of PH, regardless of whether the landmarks are intended for use in combination with the lazy witness complex or simply as a subset of the dataset to apply PH to, the aim of a landmark selection method should be to represent the underlying topology of the dataset. Since outliers can artificially introduce topological features such as loops, we in particular require outlier-robust landmark selection techniques. Based on these observations we formulate the following goals for landmark selection methods intended for PH:
\begin{enumerate}
\item Good representation of the underlying topology of the dataset, even at low sampling densities with small variance of the results between different landmark realisations.
\item Robustness to outliers, ideally including a measure for how much we consider a specific point to be an outlier.
\end{enumerate}
We now introduce novel landmark selection methods to achieve these goals: \nameref{subsec:PHOutliers} (PH landmarks) and \nameref{subsec:kminusminus}.  We design PH landmarks specifically for the application of PH, while \nameref{subsec:kminusminus} are based on a variant of the $k$--means algorithm, whose outlier-robust properties make it a promising candidate to overcome the downsides of both the random and the  maxmin landmark selection.

\subsection{Mathematical motivation}\label{sec:Motivation}

The use of our notion of local PH is inspired by the Mayer-Vietoris sequence, which can enable computation of the homology of a space $X$ by considering subspaces, whose homology is easier to compute. 
 The Mayer-Vietoris sequence for a topological space $X$ is a special type of long exact sequence. By long exact sequence we mean a sequence of abelian groups $A_i$ and group homomorphisms $\Phi_i$, $i \in \Z$,
 of the form
 
 \begin{equation}
 \dots \rightarrow \mathcal{A}_i\xrightarrow{\Phi_i} \mathcal{A}_{i+1}\xrightarrow{\Phi_{i+1}}\mathcal{A}_{i+2}\xrightarrow{\Phi_{i+2}} \dots
 \end{equation}
  such that $\Ima {\Phi}_i = \ker {\Phi}_{i+1}$. For more definitions and background, see for example  \citep{munkres1984}).

\begin{theorem}[Mayer-Vietoris sequence]
Let $X$ be a simplicial complex with subcomplexes $A,B\subset X$ such that $X = A \cup B$. Then there exists an exact sequence \\
\begin{equation}\label{eq:MVSeq}
 \begin{split}
\dots \rightarrow H_n(A \cap B)\xrightarrow{\Phi_*} H_n(A) \oplus H_n(B)\xrightarrow{\Psi_*} H_n(X)\xrightarrow{\partial_*} H_{n-1}(A \cap B)\rightarrow \dots \\ \rightarrow H_0(X)\rightarrow0,
 \end{split}
\end{equation}
where $H_n(\cdot)$ denotes the $n$-th homology group and $\Phi_*$, $\Psi_*$, and $\partial_*$ are homomorphisms.
The sequence is called the \emph{Mayer-Vietoris sequence}.
\end{theorem}

For a proof see, for example, \citep{munkres1984}, page 142.
We can use the Mayer-Vietoris sequence to connect the homology of a simplicial complex $X$ to the local homology around a vertex $\hat{x} \in X$. We do this via two simplicial subcomplexes that can be defined around the vertex $\hat{x}  \in X$: the link of the vertex $\hat{x} \in X$ and the closed star of the vertex $\hat{x}  \in X$. 

\begin{definition}[Closed star of a vertex]
Let $X$ be a simplicial complex and $\hat{x}  \in X$ a vertex. Then the closed star of $ \hat{x}$ in $X$, denoted by $\overline{\text{St}} \  \hat{x}$, is the subcomplex of $X$ that contains all the simplices which have $\hat{x}$ as one of their vertices.
\end{definition}

\begin{definition}[Link of a vertex]
Let $X$ be a simplicial complex and $\hat{x}  \in X$ a vertex. Then the link $\text{Lk} \  \hat{x}$ is the union of all simplices of $X$ lying in $\overline{\text{St}} \  \hat{x}$ that are disjoint from $\hat{x}$.
\end{definition}

\begin{figure}[htp]
\centering
\includegraphics[width=\textwidth]{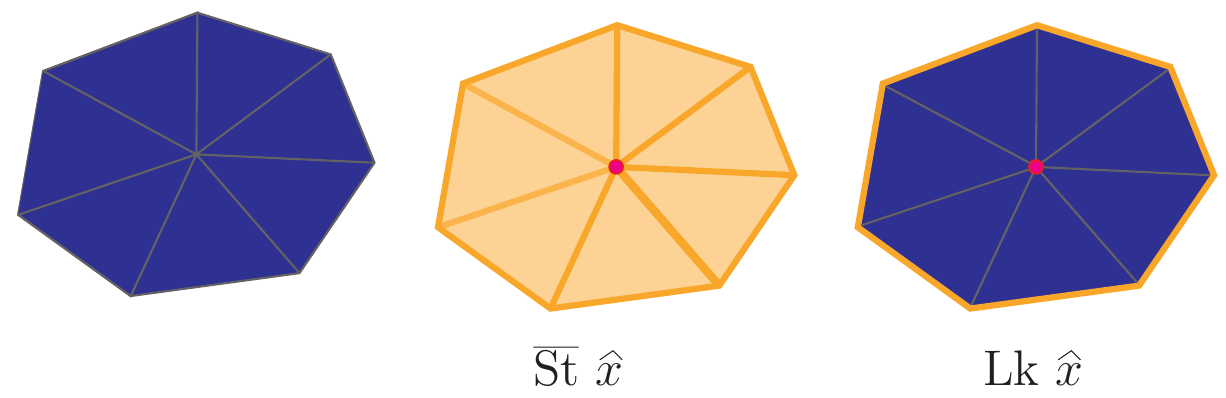}
\caption[Examples of a simplicial complex, the closed star of a vertex $\hat{x}$, and the link of a vertex $\hat{x}$.]{\label{StarLink} Examples of a simplicial complex, the closed star of a vertex $\hat{x}$, and the link of a vertex $\hat{x}$. We show the vertex $\hat{x}$ in red and highlight the closed star and the link in yellow.}
\end{figure}

\noindent We show a simplicial complex, the closed star of a vertex and the link of a vertex in Fig.~\ref{StarLink}.
Denoting the simplicial complex of all simplices in $X$ that are disjoint from $\hat{x}$ as $X\setminus \hat{x}$, we make the following observations for $\hat{x}  \in X$:
\begin{remark} \
\begin{enumerate}
\item $X = (X\setminus \hat{x})  \cup \overline{\text{St}} \  \hat{x}$.
\item $\text{Lk} \ \hat{x} = (X\setminus \hat{x} )\cap \overline{\text{St}} \  \hat{x}$.
\end{enumerate}
\end{remark}

\noindent Based on these definitions and observations we can now consider the following Mayer-Vietoris sequence:
\begin{equation}\label{eq:MVSeqStarLink}
\begin{split}
\dots \rightarrow H_n(\text{Lk} \ \hat{x})\xrightarrow{\Phi} H_n(X\setminus \hat{x}) \oplus H_n(\overline{\text{St}} \  \hat{x})\xrightarrow{\Psi} H_n(X)\xrightarrow{\partial} H_{n-1}(\text{Lk} \ \hat{x})\rightarrow \dots \\ \rightarrow H_0(X)\rightarrow0.
\end{split}
\end{equation}
Now, $ \overline{\text{St}} \  \hat{x} $ is contractible: every simplex that contains $\hat{x}$ is contractible and the intersection of simplices in $ \overline{\text{St}} \  \hat{x} $ is either empty or a simplex that contains $\hat{x}$ and is hence also contractible.
Thus $H_n( \overline{\text{St}} \  \hat{x}) = 0$ for $n>0$. We observe that if we can ensure that  $H_n(\text{Lk} \ \hat{x}) = H_{n-1}(\text{Lk} \ \hat{x}) = 0$, for $n>0$ we obtain:
\begin{equation}\label{eq:MVSeqStarLinkII}
0\xrightarrow{} H_n(X\setminus \hat{x}) \xrightarrow{\Psi} H_n(X)\xrightarrow{} 0,
\end{equation}
which gives us an isomorphism $\Psi$ between $H_n(X\setminus \hat{x})$ and $H_n(X)$. 

The Mayer-Vietoris sequence given by Eq.~\ref{eq:MVSeqStarLink} connects the homology of the large simplicial complexes $X$ and $X\setminus \hat{x}$, which are both global and expensive to compute, to the  homology of the Link $\text{Lk} \ \hat{x}$, which a purely local computation and therefore easy to compute. 

In practice, we are however working with point cloud data $D$. We therefore consider the PH of the point cloud $D$, rather than the homology, and apply the Mayer-Vietoris sequence~\ref{eq:MVSeqStarLink} to the simplicial complexes that we obtain from the data by constructing a filtration. We consider a simplicial complex $X$ in this filtration. A data point $y \in D$ is now a vertex $\hat{y}$ in $X$. 
Instead of ${H_n(\text{Lk} \ \hat{y}) =  0}$, we need to quantify the failure of the isomorphism ${PH_n(D\setminus y) \xrightarrow{\Psi} PH_n(D)}$ via a score based on the $PH_n(\text{Lk} \ y)$, which we will 
define in Section \ref{subsec:PHOutliers}. 
To make computation easier, instead of looking just at the link of $\hat{y}$ in the simplicial complex $X$, we extend the link to a $\delta$-neighbourhood of $\hat{y}$ in $X$, which we define to be the collection of simplices whose vertices are within a distance of at most $\delta$ from $y$ in $D$: 

\begin{definition}[$\delta$-link of a data point $y$]
Let $X$ be a simplicial complex in a metric space constructed from a dataset $D$, $\hat{y}  \in X$ a vertex, $\delta>0$ a distance and $\delta(\hat{y})$ a $\delta$-neighbourhood of $\hat{y}$ in $X$. Then the $\delta$-link $\text{Lk}^\delta \ y$ is the union of all simplices in $X$ that are disjoint from $\hat{y}$ and contained in $\delta(\hat{y})$ .
\end{definition}

\noindent Building on the $\delta$-Link of a data point we can define the $\delta$-Star of a data point:

\begin{definition}[Closed $\delta$-star of a data point $y$]
Let $X$ be a simplicial complex in a metric space constructed from a dataset $D$, $\hat{y}  \in X$ a vertex, $\delta>0$ a distance and $\delta(\hat{y})$ a $\delta$-neighbourhood of $\hat{y}$ in $X$. Then the closed $\delta$-star $\overline{\text{St}}^\delta  \ \hat{y}$ is the union of the $\delta$-link $\text{Lk}^\delta \ y$ with the vertex $\hat{y}$ and all simplices $[\hat{y},\sigma]$ where $\sigma \in X$ and $\sigma$ is fully contained in $\delta(\hat{y})$.
\end{definition}

\begin{figure}[htp]
\centering
\subcaptionbox{$\delta$-neighbourhood of $y$ in $D$.}{\includegraphics[width=.3\textwidth]{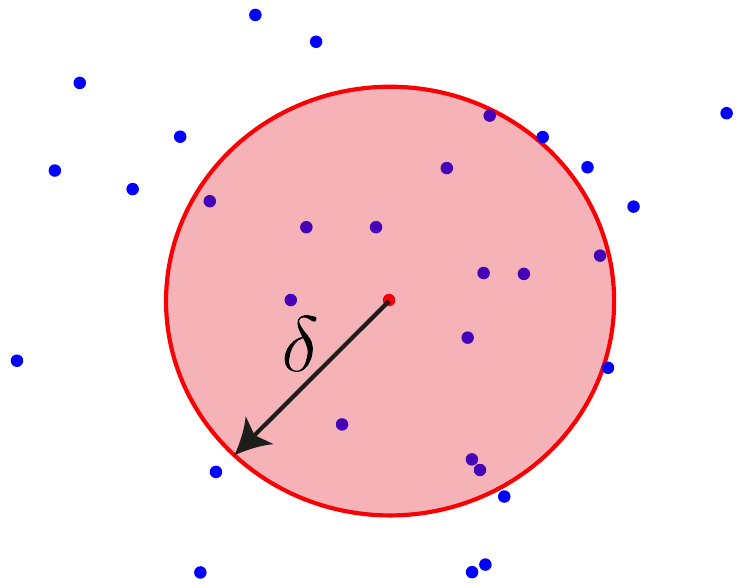}}%
\hspace{0.1cm}
\subcaptionbox{$\text{Lk}^\delta \ y$ at scale $\epsilon$.}{\includegraphics[width=.3\textwidth]{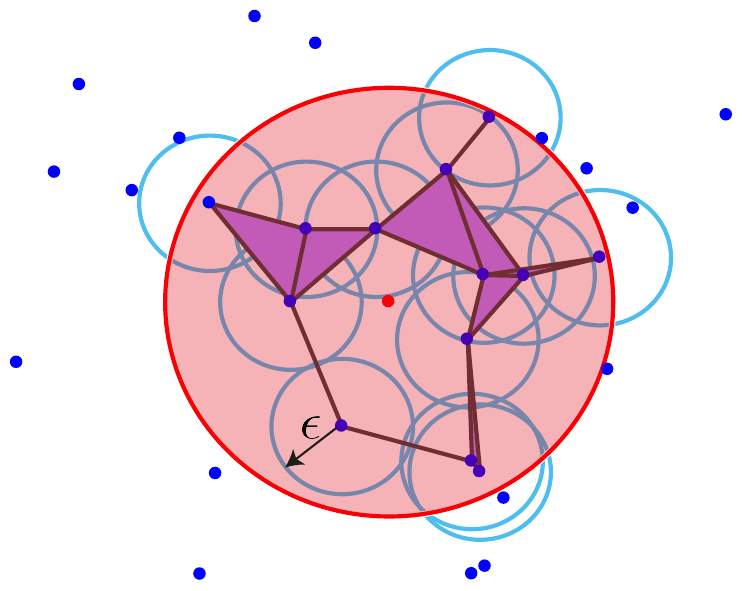}}%
\hspace{0.1cm}
\subcaptionbox{$\overline{\text{St}}^\delta \ y$ at scale $\epsilon$.}{\includegraphics[width=.3\textwidth]{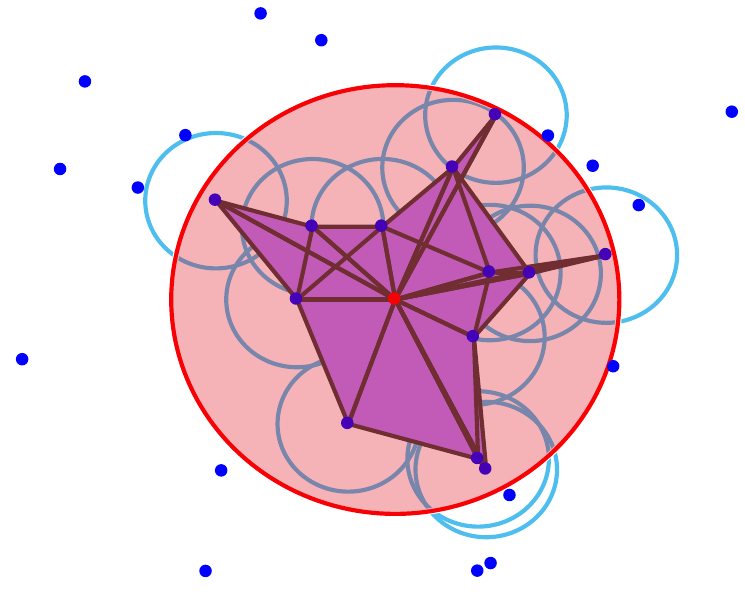}}%
\caption[Examples of a data point $y$ and its $\delta$-neighbourhood in a point cloud, the $\delta$-link of $y$ and  the closed $\delta$-star of $y$. ]{\label{DeltaStarLink} Examples of a data point $y$ and its $\delta$-neighbourhood in a point cloud, the $\delta$-link of $y$ and the closed $\delta$-star of $y$. We show the data point $y$ in red and its $\delta$-neighbourhood in light red highlighting the points within $\delta$ in blue.
We use the data points within the $\delta$-neighbourhood to build a Vietoris--Rips complex for a set filtration value $\epsilon > 0$, which represents the simplicial complex $X$ in the definitions of the $\delta$-link and the closed $\delta$-star. We only show the subcomplexes of the Vietoris--Rips complex (or their extensions in the case of the closed $\delta$-star) that are relevant to the illustrated definitions.}
\end{figure}

We note that the closed $\delta$-star of a data point is always contractible by construction.
We show an example of a data point, its $\delta$-link and its closed $\delta$-star in Fig.~\ref{DeltaStarLink}.
From now on, we refer to computing the PH of the $\delta$-Link of a point in a dataset as computing the local PH of a data point. Our notion of local PH of a data point is motivated by the Mayer-Vietoris sequence and we use the property given by Eq.~\ref{eq:MVSeqStarLinkII} to define a new landmark selection method, in which we select points $ y$ in a dataset $D$ which are `closest' to giving us the desired isomorphism between $PH_n(D\setminus y)$ and $PH_n(D)$ as landmarks. Note that a similar computation for local PH is performed in \cite{Wheeler2021}.

\subsection{Persistent homology landmarks}\label{subsec:PHOutliers}

Following our local PH computations, we can now define the score $|PH_n(\text{Lk}^\delta \ y)|$ by which we measure the failure of the isomorphism ${PH_n(D\setminus y) \xrightarrow{\Psi} PH_n(D)}$. 

\begin{definition}[$|PH_n(\text{Lk}^\delta \ y)|$]\label{def:PH_out}
Let $D$ be a point cloud, $y  \in D$ a data point, $d:D\times D \rightarrow \R$ a distance function, $\Delta_{y} = \{ \tilde{y} \in D\setminus{\{y \}} \ | \ d(\tilde{y},y) \leq \delta \}$, $n = 0,1,2$ and $\mathcal{B}_n(\text{Lk}^\delta \ y) = \{ [\eta_i ,\zeta_i) \}_{i = 1}^{I(n)}$ the $n $-dimensional barcode of the Vietoris--Rips filtration performed on $\Delta_{y}$ excluding infinitely persisting features. Then
$$|PH_n(\text{Lk}^\delta \ y)| = \max\limits_{n}{\max \limits_{i = 1,\dots,I(n)} }\{  \zeta_i-\eta_i \}.$$   
\end{definition}

We can reformulate our observation from Eq.~\ref{eq:MVSeqStarLinkII} for point cloud data $D$ and PH in the following way: if for a data point $y  \in D$ and a neighbourhood radius $\delta>0$ we can ensure that  $|PH_n(\text{Lk}^\delta \ y)| = |PH_{n-1}(\text{Lk}^\delta \ y)| = 0$, then for $n>0$ we obtain:
\begin{equation}\label{eq:MVSeqStarLinkIII}
0\xrightarrow{} PH_n(D\setminus y) \xrightarrow{\Psi} PH_n(D)\xrightarrow{} 0,
\end{equation}
where $PH_n$ denotes the $n$-th homology groups associated with the different filtration steps. As we are working with data however, we are unlikely to achieve such strict conditions -- indeed, we would achieve these under trivial conditions. For $|PH_n(\text{Lk}^\delta \ y)| = |PH_{n-1}(\text{Lk}^\delta \ y)| \approx \text{small}$, however, we obtain something close to an isomorphism between $PH_n(D\setminus y) $ and $PH_n(D)$. 

The larger $|PH_n(\text{Lk}^\delta \ y)|$ of a point is, the further away we are from an isomorphism between $PH_n(D\setminus y) $ and $PH_n(D)$. Consequently, the inclusion or exclusion of the point changes the PH of the dataset more than for a point where $|PH_n(\text{Lk}^\delta \ y)|$ is small. 
Under the assumption that each point has at least two neighbours within distance $\delta$, the choice of $\delta$ determines by how much we allow the Bottleneck distance (for a definition see, for example, \citep{Otter2017}) between the persistence diagrams of the filtration on the point cloud containing $y$ and a point cloud not containing $y$ to differ. We can therefore think of it as a resolution parameter.
There are now two possible strategies for landmark selection:
\begin{itemize}
\item \emph{PH landmarks I: representative landmarks.} Landmarks are points with small $|PH_n(\text{Lk}^\delta \ y)|$ values. Their inclusion or exclusion in the full dataset does not alter the PH of the full dataset dramatically. Hence, the landmarks represent the data well and they do not introduce accidental topological features that can, for example, be caused by outlier points in the landmark set.
\item \emph{PH landmarks II: vital landmarks.} Landmarks are points with large $|PH_n(\text{Lk}^\delta \ y)|$ values. For such points we are far away from an isomorphism between $PH_n(D\setminus \hat{y}) $ and $PH_n(D)$. Hence their exclusion from the dataset would change the PH dramatically.
\end{itemize}

\noindent In practice, it is not immediately clear which strategy works better. Here, we therefore consider both.
We compute $|PH_n(\text{Lk}^\delta \ y)|$ for dimensions $n=0,1,2$. 
Motivated by results from our computations, we define two versions for $|PH_n(\text{Lk}^\delta \ y)|$ which we use in practice. To avoid confusion with Definition \ref{def:PH_out}, we refer to these values obtained in practical computations as the \emph{PH outlierness} of a point $y \in D$:
\begin{equation}\label{eq:out_all_dim}
out_{\text{PH}}^{0,1,2}(y) = \max{ \{|\mathcal{B}_0(\text{Lk}^\delta \ y)|,|\mathcal{B}_1(\text{Lk}^\delta \ y)|,|\mathcal{B}_2(\text{Lk}^\delta \ y)|\}},
\end{equation}
where $|\mathcal{B}_n(\text{Lk}^\delta \ y)|$ is the length of the longest finite interval in the PH barcode for dimension $n$. We also refer to $out_{\text{PH}}^{0,1,2}(y)$ as outlierness over all (computed) dimensions.
We observe that $out_{\text{PH}}^{0,1,2}(y)$ is usually determined by dimension 0, where we find the longest non-infinitely persisting features in our datasets (see Subsection \ref{sec:caseStudy}). To avoid this, we also use a version for PH outlierness which is restricted to dimension 1 in the PH calculation:
\begin{equation}\label{eq:out_dim1}
out_{\text{PH}}^{1}(y) = |\mathcal{B}_1(\text{Lk}^\delta \ y)|.
\end{equation}
We also refer to $out_{\text{PH}}^{1}(y)$ as dimension 1 outlierness. In situations where either of the two definitions can be used, we denote PH outlierness as $out_{\text{PH}}(y)$.

To avoid choosing points as landmarks that are very far away from other data points, we determine points $S = \{s_1,\dots,s_o\}$ with fewer than two neighbours within their $\delta$-neighbourhood to be \emph{super outliers}. We include super outliers into the landmark set only once all other points have been chosen as landmarks. Note that the resolution parameter $\delta$ strongly influences the number of super outliers. As long as we have enough points in the dataset that are not considered to be super outliers, we choose our landmarks to be the points $L = \{ l_1, l_2,\dots,l_m\} \subset D\setminus S$ such that ${out_{\text{PH}}(l_i) \leq out_{\text{PH}}(y)}$ for all $y \in D\setminus \{L \cup S\}$ and $i = 1, \dots, m$ for PH landmarks I (representative landmarks) and ${out_{\text{PH}}(l_i) \geq out_{\text{PH}}(y)}$ for all $y \in D\setminus \{L \cup S\}$ and $i = 1, \dots, m$ for PH landmarks II (vital landmarks). 
We show the pseudocode for PH landmarks I in Algorithm~\ref{alg:PHLandmarks}, an algorithm for PH landmarks II can be formulated accordingly.
\begin{algorithm}
\caption{The PH landmark algorithm}
\label{alg:PHLandmarks}
\begin{algorithmic} 
\REQUIRE Data points $D = \{y_1,\dots,y_N\}$, \\ a distance function $d:D\times D \rightarrow \R$ \\ number of landmarks $m$,\\ local neighbourhood radius $\delta>0$.
    \ENSURE A set of $m$ PH landmarks $L = \{l_1,\dots,l_m\}$, a set of $o$ super outliers ${S = \{s_1,\dots,s_o\}}$.
\FORALL{$y \in D$}  
\STATE Find $\Delta_{y} = \{ \tilde{y} \in D\setminus{\{y \}} \ | \ d(\tilde{y},y) \leq \delta \}$ 
\IF{$|\Delta_{y}| >1$}
\STATE Compute Vietoris--Rips filtration for $\Delta_{y}$ for $n = 0,1, 2$.
\STATE Compute $out_{\text{PH}}(y)$
\ELSE
\STATE $S \leftarrow S \cup \{y\}$
\ENDIF
\ENDFOR
\STATE Re-order the points in $D\setminus S$ such that  $out_{\text{PH}}(y_1) \leq out_{\text{PH}}(y_2) \leq \dots \leq out_{\text{PH}}(y_{N-o})$
\STATE $L \leftarrow \{y_1,\dots,y_{\min\{m,{N-o}\}}\}$
\IF{$N-o<m$}
\STATE  $L \leftarrow L \cup \{s_1,\dots, s_{m-N+o}\}$
\ENDIF
\end{algorithmic}
\end{algorithm}

When using $out_{\text{PH}}^{1}(y)$ in Algorithm~\ref{alg:PHLandmarks}, it is important to note, that one can obtain data points $y$ that are not super outliers with $out_{\text{PH}}^{1}(y) = 0$. To avoid that the order of inclusion of such points in the landmark set is determined by a possibly systematic and non-random ordering of the points in the original dataset, which could, for example, favour noise points to be added before signal points or vice versa, we ensure that all points with $out_{\text{PH}}^{1}(y) = 0$ are randomly permuted in the ordering of the PH outlierness scores.

\subsection{$k--$ landmarks}\label{subsec:kminusminus}

The $k$--means$--$ algorithm was developed by 
\citep{Chawla2013} to overcome the extreme sensitivity of the $k$--means algorithm to outliers. The authors formulate their approach as a generalisation of the $k$--means algorithm: for an input dataset $D = \{y_1,\dots,y_N\}$ the algorithm provides a set of $k$ cluster centres $\hat{L} = \{\hat{l}_1,\dots,\hat{l}_k\}$ and a set of $j$ outliers $O = \{o_1,\dots,o_j \}$, $O \subset D$. For a given distance function $d:D\times D \rightarrow \R$ and $y \in D$ the authors use the following term in their algorithm:
\begin{equation}
c(y, \hat{L}) \coloneqq \text{arg}\min_{\hat{l} \in \hat{L}}{ d(y,\hat{l})}.
\end{equation}

\noindent We show the pseudocode in Algorithm~\ref{alg:kminusminus} in Section~\ref{sec:Pseudocodekminusminus}.  For our application of the algorithm to landmark selection, we further define:
\begin{equation}
\tilde{c}(D,\hat{l}) \coloneqq \text{arg}\min_{y \in D}{ d(y,\hat{l})}.
\end{equation}
\noindent  We show our modified version of the $k$--means$--$ algorithm for landmark selection in Algorithm~\ref{alg:kminusminusExt}.
\begin{algorithm}
\caption[The $k$--means$--$ algorithm modified for landmark selection]{The $k$--means$--$ algorithm~\citep{Chawla2013} modified for landmark selection}
\label{alg:kminusminusExt}
\begin{algorithmic} 
    \REQUIRE Data points $D = \{y_1,\dots,y_N\}$, a distance function $d:D\times D \rightarrow \R$, number of clusters $k$ and number of outliers $j$.
    \ENSURE A set of $k$ cluster centers $L = \{l_1,\dots,l_k\}$, \ $  L \subset D$, \\ a set of $j$ outliers $O = \{o_1,\dots,o_j \}$, $O \subset D$.
    \STATE $\hat{L}_0 \leftarrow \{k \text{ random points of } D\}$
      \STATE $e_0 = -1$  \color{black}
    \STATE $i \leftarrow 1$
    \WHILE{ ($\text{continuation\_criterion} >  10^{-4}$  \AND $i < 100$ )\color{black}  }  
            \FORALL{$y \in D$}  \STATE compute $d(y,\hat{l}_{i-1})$ \ENDFOR
            \STATE Re-order the points in $D$ such that $d(y_1, \hat{L}_{i-1}) \geq d(y_2, \hat{L}_{i-1}) \geq \dots \geq d(y_N, \hat{L}_{i-1}) $
            \STATE $O_i \leftarrow \{y_1,\dots,y_k\}$
            \STATE $D_i \leftarrow D\setminus O_i = \{y_{k+1},\dots,y_N\}$
            \FOR{r = 1 \TO k}
            	\STATE {$P_r \leftarrow \{y \in D_i \ | \ c(y, \hat{L}_{i-1}) = \hat{l}_{i-1,r} \}$}
		\STATE {$\hat{l}_{i,r} \leftarrow \text{mean}(P_r)$}
            \ENDFOR
           \STATE $\hat{L}_i \leftarrow \{\hat{l}_{i,1}, \dots, \hat{l}_{i,k}\}$
          
   	 \FOR {$y \in D_i$}
    		\STATE{compute $d(y,\hat{L}_i)$}
    	\ENDFOR
	\STATE $e_i \leftarrow \sum_{y\in D_i}{d(y,\hat{L}_i)^2}$
	\STATE $\text{continuation\_criterion } \leftarrow |e_i-e_{i-1}|$
	\color{black}
           \STATE $i \leftarrow i + 1$
    \ENDWHILE
     
    \STATE $L \leftarrow \emptyset$
    \STATE $D_L \leftarrow D_{i-1}$ 
    \WHILE {$|L| < k$}
    
    	\FOR{$\hat{l} \in \hat{L}$}
    	\STATE $ m_{\hat{l}}\leftarrow \min_{y \in D_{L}}{d(y,\hat{l})}$
    	\ENDFOR  
	\STATE Re-order $\{\hat{l}_{1}, \dots, \hat{l}_{\hat{k}}\}$ such that $m_{l_{1}}\leq m_{l_{2}}\leq \dots \leq m_{l_{\hat{k}}}$
	\STATE $s \leftarrow 1$
	\REPEAT 
	\STATE$L \leftarrow L \cup \{ y_s\}$, where $y_s = \tilde{c}(D_{L},\hat{l}_{s})$
	\STATE $s \leftarrow s + 1$
	\UNTIL{$y_s = y_t$ for some $t<s$}
	\STATE $L \leftarrow L\setminus \{y_s\}$
	\STATE $D_L \leftarrow D_L \setminus L$
	\STATE $\hat{L} \leftarrow \hat{L}\setminus \{ \hat{l}_1,\dots ,\hat{l}_{s-1}\}$
     \ENDWHILE
    
    \color{black}

\end{algorithmic}
\end{algorithm}

\subsection{Implementation}\label{subsec:LandmarkImplementation}

We implement PH landmark selection method in {\sc Matlab} using {\sc Ripser}~\citep{bauer2017ripser,Bauer2021} for the computation of the local Vietoris--Rips complexes. We also implement the $k--$ landmarks algorithm in {\sc Matlab}. For the calculation of  maxmin landmarks, random landmarks, and dense core subsets we use the inbuilt functions in the {\sc javaPlex} package~\citep{Javaplex}.

\section{Datasets}\label{subsec:LandmarkData}

We introduce the datasets to which we apply the different landmark selection methods. The datasets are chosen to be simple to allow us to determine effects of the landmark selection. The datasets consist of signal points that are sampled from a topologically interesting structure -- a sphere, a Klein bottle, or a torus -- and noise points that we design to be topologically different from the signal.

\subsection{3-dimensional datasets}

\subsubsection{Sphere-cube dataset.}

For a given number of points $N$ and probability $p$ we sample signal points uniformly at random from the surface of the unit sphere with probability $p$ and noise points from the (filled) cube $[-1,1]^3\subset \R^3$ with probability $1-p$.

\subsubsection{Sphere-plane dataset}

For a given number of points $N$ we sample signal points uniformly at random from the surface of the unit sphere with probability $p$ and noise points from the $xy$-plane $[-3,3]^2 \subset \R^2$ with probability $1-p$.

\subsubsection{Sphere-line dataset}

For a given number of points $N$ we sample signal points uniformly at random from the surface of the unit sphere with probability $p$ and noise points from $(\alpha, 0, 0)$, where $\alpha \in [-50,50] \subset \R$, with probability $1-p$.

\subsubsection{Sphere-Laplace line dataset}

For a given number of points $N$ we sample signal points uniformly at random from the surface of the unit sphere with probability $p$ and we sample noise points from $(\alpha, 0, 0)$ with probability $1-p$, where $\alpha$ is sampled from $[-50,50]\subset \R$ and Laplace distributed with $\mu = 4$ and $\sigma = 0.5$. We use the Laplacian random number generator code~\citep{Laplace} to generate $\alpha$.

\subsection{4-dimensional datasets}

\subsubsection{Torus dataset}

We use the following parametrisation of the torus $\mathcal{T}$:
\begin{equation*}
(x,y,z,\omega) = (\cos(\gamma), \sin(\gamma),\cos(\varphi),\sin(\varphi)),
\end{equation*}
where $\gamma,\varphi \in (0,2\pi)$. We add noise $T_{\text{noise}}$ to the torus using the equation
\begin{equation*}
(x_{\text{noise}},y_{\text{noise}},z_{\text{noise}},\omega_{\text{noise}}) = (r*\cos(\gamma), r*\sin(\gamma), \hat{r}*\cos(\varphi), \hat{r}*\sin(\varphi)),
\end{equation*}
 where $r, \hat{r} \in (0,2)$. 

For a given number of points $N$ we sample signal points uniformly at random from $\mathcal{T}$ with probability $p$ and noise points from $\mathcal{T}_{\text{noise}}$ with probability $1-p$.

\subsubsection{Klein bottle dataset}

We use the following parametrisation of the Klein bottle $\mathcal{K}$:
\begin{align*}
x &= \cos(\gamma)*(r*\cos(\varphi)+C), \\
y &= \sin(\gamma)*(r*\cos(\varphi)+C),\\
z &= \cos(\gamma/2)*r*\sin(\varphi),\\
\omega &= \sin(\gamma/2)*\sin(\varphi),
\end{align*}
where $\gamma,\varphi \in (0,2\pi)$, $r = 3$ and $C=2$. We define noise $\mathcal{K}_{\text{noise}}$ for the Klein bottle using the equations:
\begin{align*}
x_{\text{noise}} &= \cos(\gamma)*(r_{\text{noise}}*\cos(\varphi)+C_{\text{noise}}), \\
y_{\text{noise}} &= \sin(\gamma)*(r_{\text{noise}}*\cos(\varphi)+C_{\text{noise}}),\\
z_{\text{noise}} &= \cos(\gamma/2)*r_{\text{noise}}*\sin(\varphi),\\
\omega_{\text{noise}} &= \sin(\gamma/2)*\sin(\varphi),
\end{align*}
where $r_{\text{noise}}$ is sampled uniformly from the interval $[2,4]\subset \R$ and $C_{\text{noise}}$ is sampled uniformly from the interval $[1,3]\subset \R$. We use and adapt the code from~\citep{Otter2017}. For a given number of points $N$ we sample signal points uniformly at random from $\mathcal{K}$ with probability $p$ and noise points from $\mathcal{K}_{\text{noise}}$ with probability $1-p$.

\section{Results}\label{subsec:LandmarkResults}

We present our results for the proposed landmark selection methods, i.e. PH landmarks and $k--$ landmarks. 
We first study PH landmark selection in detail on the sphere-cube dataset with $p = 0.6$, then proceed to showing our results in comparison to the current standard methods on our various datasets, i.e. landmark selection via the maxmin algorithm and random landmark selection. We then compare our methods to dense core subsets and finally investigate the influence of the $\delta$ parameter on the performance of PH landmarks. For all landmark selection techniques we use the Euclidean distance as distance function. All of our datasets consist of 3000 points. 

\subsection{PH landmarks case study on the sphere-cube dataset with $p = 0.6$}\label{sec:caseStudy}

\begin{figure}[htp]
\centering
\includegraphics[width=\textwidth]{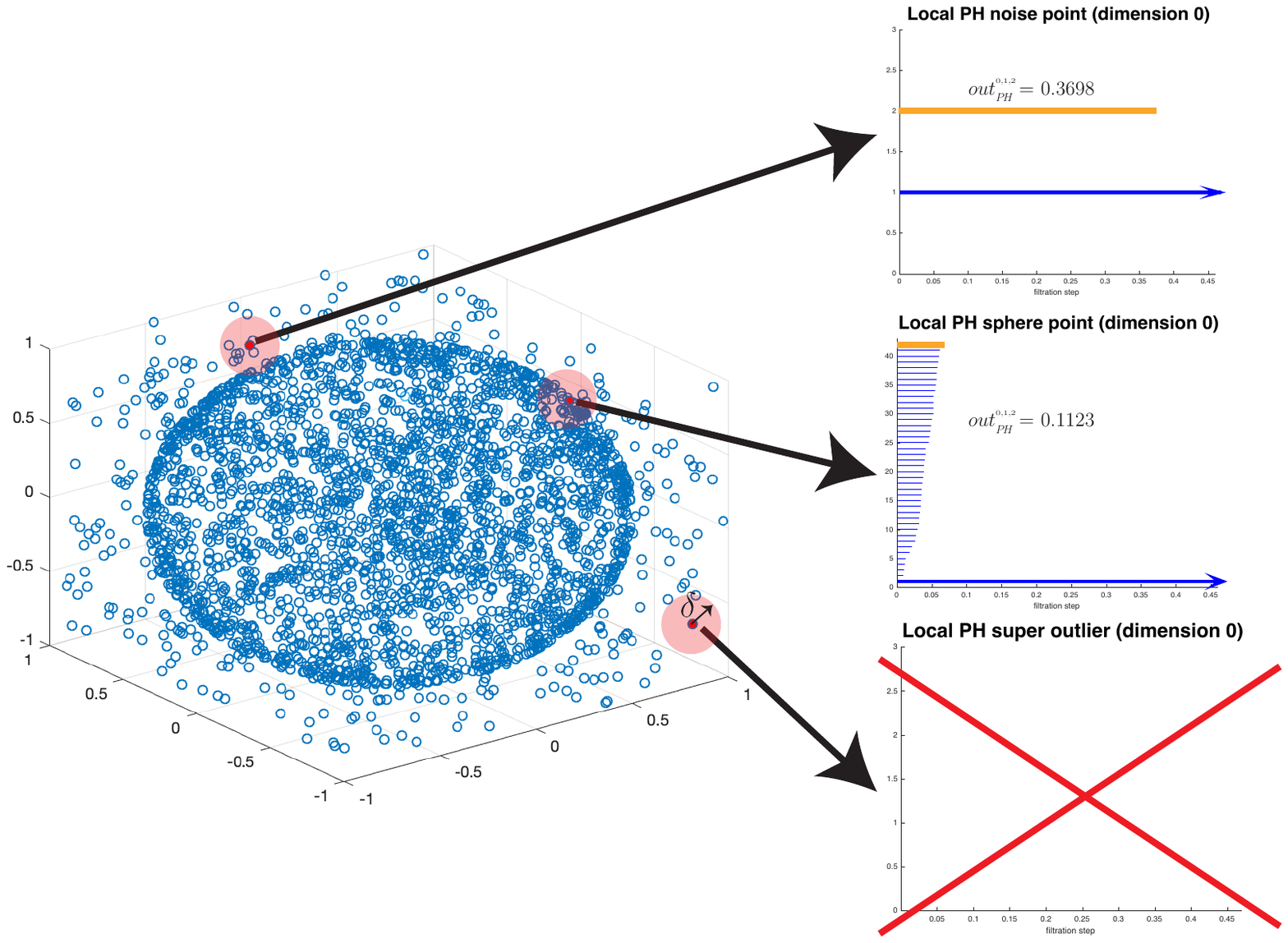} 
\caption[Schematic illustration of three different types of points found on the sphere-cube dataset]{\label{SphereCubeMethod} Schematic illustration of three different types of points $y$ found in the sphere-cube dataset, $p = 0.6$, and their local dimension $0$ barcodes for $\delta = 0.2$: noise point, sphere point, and super outlier. The $\delta$-neighbourhoods are shown in light red balls around the corresponding data points. We calculate the PH outlierness $out_{\text{PH}}^{0,1,2}(y)$ for every point based on its local Vietoris--Rips barcode and ignore super outliers. We highlight the bars in the barcodes that are used to determine $out_{\text{PH}}^{0,1,2}(y)$ in yellow.}
\end{figure}

\begin{figure}[htp]
\centering
\includegraphics[width=0.5\textwidth]{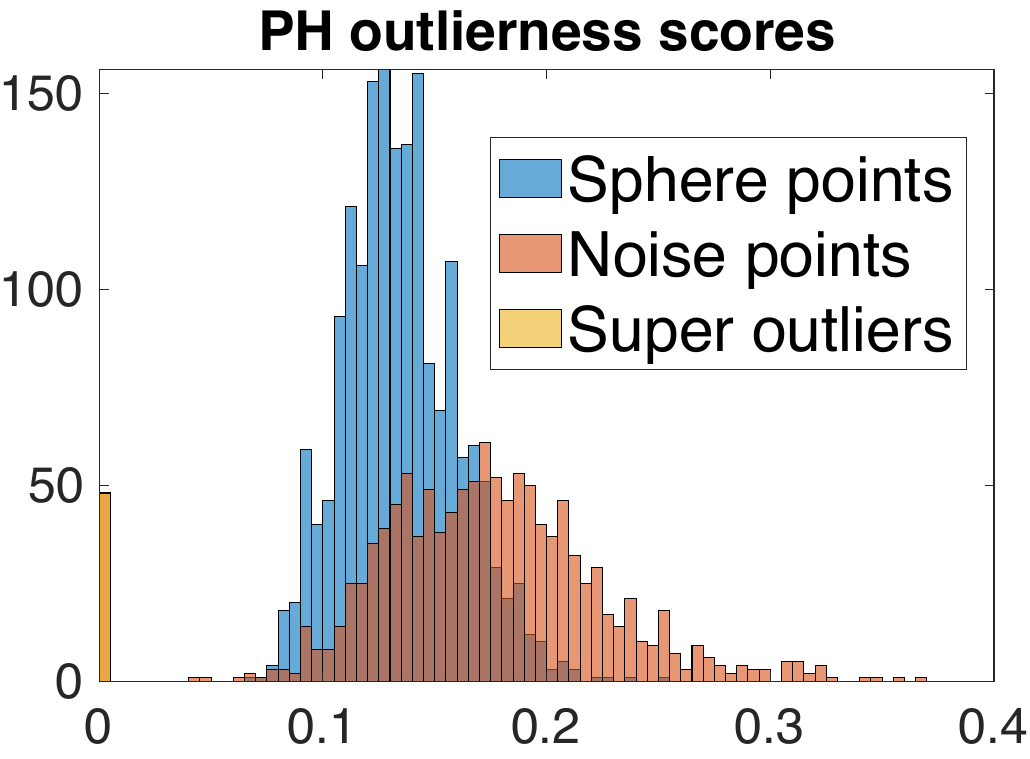}
\caption{\label{SphereCubeHistogram} Histograms of the PH outlierness $out_{\text{PH}}^{0,1,2}(y)$ values obtained on the sphere-cube dataset, $p = 0.6$, from local PH with $\delta = 0.2$. The horizontal axis represents the outlierness scores, the vertical axis shows the number of points.}
\end{figure}

We apply PH landmark selection with $\delta = 0.2$ to the sphere-cube dataset where a data point has a probability of 0.6 to be located on the surface of the unit sphere and 0.4 to be located in the unit cube. As for all of our datasets, we find that the PH outlierness values $out_{\text{PH}}^{0,1,2}(y)$, as defined in Eq.~\ref{eq:out_all_dim}, are determined by the maximal non-infinite bar in the dimension $0$ barcode. We show example barcodes for dimension 0 for a noise point, a sphere point, and a super outlier in Fig.~\ref{SphereCubeMethod}. In general, we expect a noise point to be located in a sparser region of this dataset and thus to either be classified as a super outlier, or to exhibit a barcode with a small number of long bars and very few, or no, short bars. For a sphere point, we expect the dimension 0 barcode to have many short bars and occasionally some longer bars caused by noise points that lie within the $\delta$-neighbourhood. Note that for the examples in Fig.~\ref{SphereCubeMethod}, we choose the noise point with the highest outlierness score and the sphere point with the lowest outlierness score in the dataset to illustrate model cases for the method.
We find that for these example points, the dimension 0 barcodes behave as expected. To explore whether the outlierness scores reflect the properties of the different types of points as expected, we consider histograms of the outlierness scores $out_{\text{PH}}^{0,1,2}(y)$ of all data points in Fig.~\ref{SphereCubeHistogram}.
We find that we have 48 super outliers in the dataset. For the noise and sphere points, we can see that the outlier scores are distributed differently and that by including points with low outlier scores as landmarks first, we should preferentially obtain sphere points rather than noise points. Landmarks chosen in this way correspond to representative landmarks (PH landmarks I) described in Subsection~ \ref{subsec:PHOutliers}.

\begin{figure}[htp]
\centering
\includegraphics[width=0.5\textwidth]{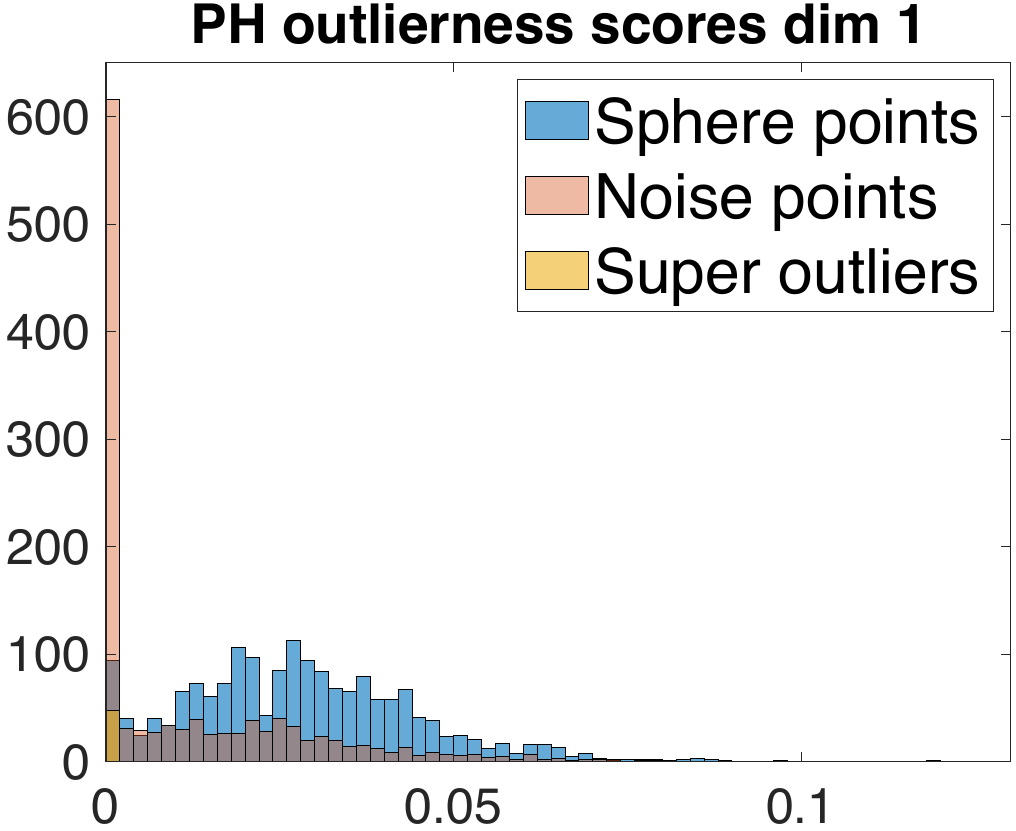}
\caption{\label{SphereCubeHistogramII} Histograms of the PH outlierness $out_{\text{PH}}^1(y)$ values obtained on the sphere-cube dataset, $p = 0.6$, from local PH with $\delta = 0.2$, considering only features in dimension 1. The horizontal axis represents the outlierness scores, the vertical axis shows the number of points.}
\end{figure}

As our landmark selection approach is motivated by Eq.~\ref{eq:MVSeqStarLinkIII} which holds for PH in dimensions $n>0$, it is not immediately clear that representative landmarks are preferable to vital landmarks both for $n = 0$ and dimensions $n>0$. We examine a variant of the method where we restrict ourselves to local PH in dimension~1. In this case dimension~1 PH outlierness values correspond to the persistence of the most persistent feature in the local dimension 1 barcode (see Eq.~\ref{eq:out_dim1}). We show a histogram of the distribution of the dimension 1 PH outlierness scores in Fig.~\ref{SphereCubeHistogramII}.
We observe that the clearest difference between the sphere points and the noise points is the fact that a large proportion of noise points has $out_{\text{PH}}^1(y) = 0$, while a clear majority of the sphere points has $out_{\text{PH}}^1(y) > 0$. This can again be explained by the fact that sphere points have more neighbours within their $\delta$-neighbourhood and therefore are more likely to form features in dimension~1. From these observations it seems that here it is more beneficial to use vital landmarks (PH landmarks II), described in Subsection~ \ref{subsec:PHOutliers}, for landmark selection based on dimension 1. 
For PH landmarks based on dimension 1 we thus choose points with large dimension 1 PH outlierness scores as landmarks and discard points with low dimension 1 PH outlierness scores from the dataset as outliers whose removal does not alter the PH of the dataset much.

\subsection{Comparison of persistent homology landmarks and $k--$ landmarks to standard landmark selection methods}\label{subsubsec:landmarkComparisons}

We compare our proposed methods for landmark selection, the $k--$ landmarks and the PH landmarks, to the current standard methods for landmark selection, i.e. random landmarks and  maxmin landmarks. Using the different techniques we choose $m$ landmarks from $N$ data points. This corresponds to a landmark sampling density of $\frac{m}{N}$. For
 the $k--$ landmarks, we define the number of clusters to be $k = p m$ and the number of outliers to be $j = (1-p) m$, where $p$ is the probability with which a point in the respective dataset was sampled from the signal data, i.e. the sphere, torus, or Klein bottle. We consider both the case where we choose the $k--$ cluster centres and the outliers found by the algorithm as our landmarks as well as the case where we only consider the $k--$ cluster centres to be landmarks. For the PH landmarks we choose $\delta = 0.2$ for the 3-dimensional datasets,  $\delta = 0.5$ for the torus data and $\delta = 0.6$ for the Klein bottle. 
In all our datasets, we find that when looking for the maximal persistence of a feature across all dimensions, the value of $out_{\text{PH}}^{0,1,2}(y)$ is exclusively determined by dimension 0. We therefore also include a variation of PH landmarks that only considers the local dimension 1 barcode for the calculation of $out_{\text{PH}}^1(y)$.
For the PH landmark version where we use all dimensions to determine the outlierness scores, we choose data points with $out_{\text{PH}}^{0,1,2}(y) \approx $ small as our landmarks (representative landmarks), for the version where we restrict ourselves to dimension 1, we choose points with large $out_{\text{PH}}^{1}(y)$ scores as landmarks (vital landmarks).

\begin{figure}[htp]
\centering
\includegraphics[width=.3\textwidth]{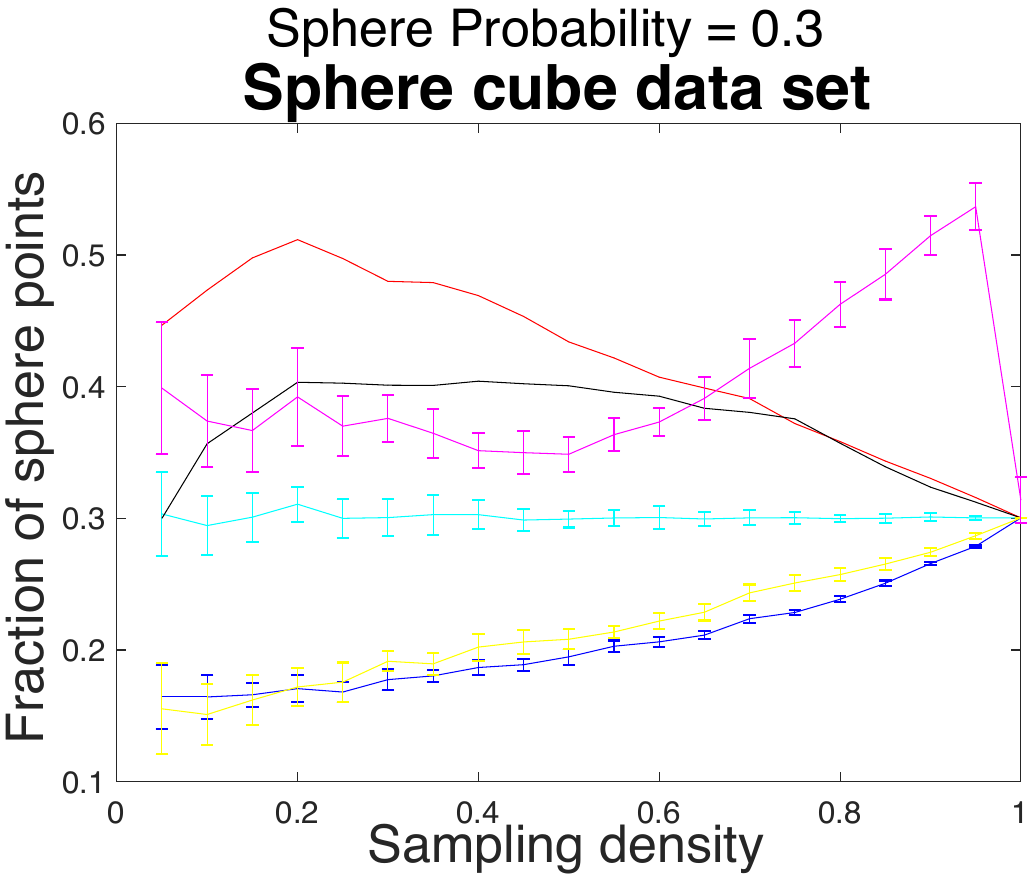}
\includegraphics[width=.3\textwidth]{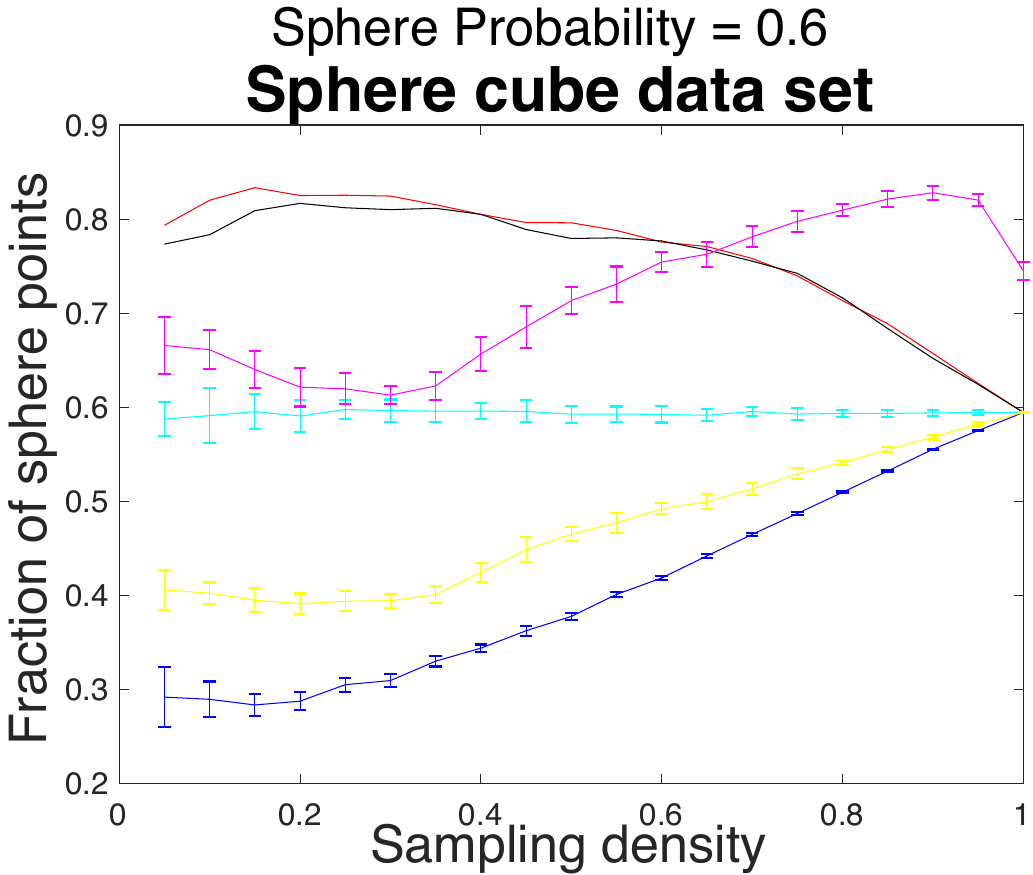}
\includegraphics[width=.3\textwidth]{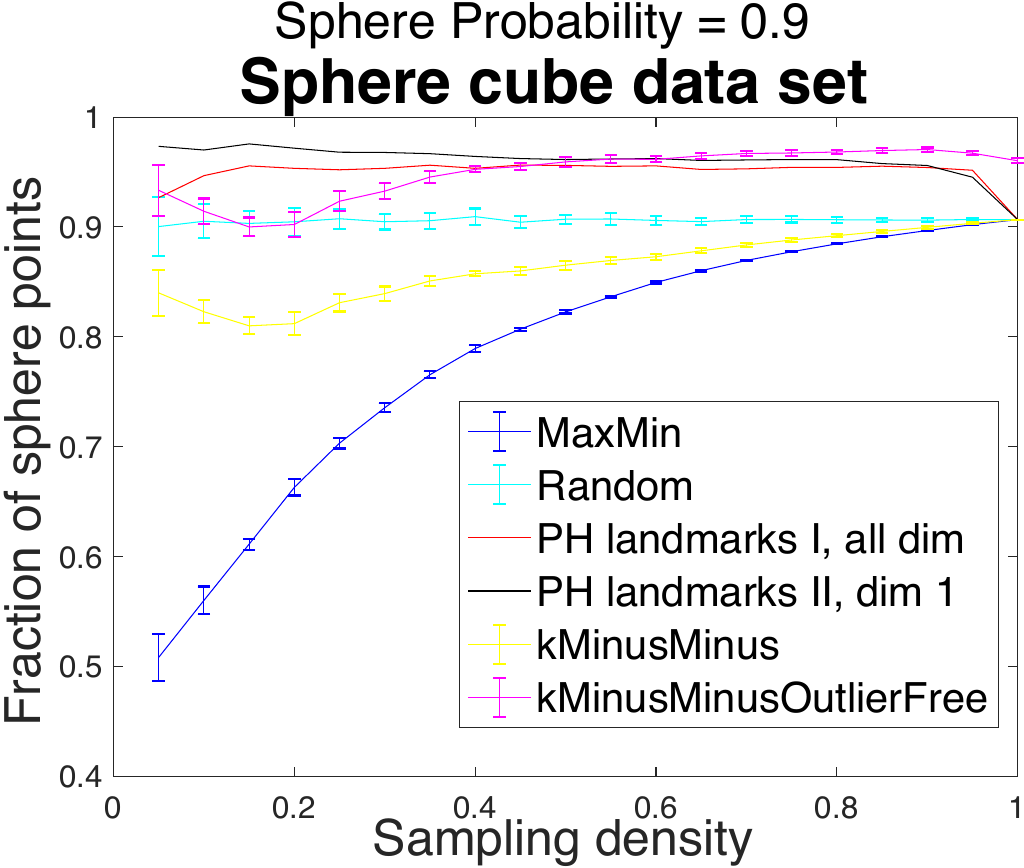}
\caption{\label{ComparisonSphereCube} Comparison of the fraction of sphere points in selected landmark points for different landmark selection techniques on the sphere-cube dataset for $\delta = 0.2$. We consider landmark selection via the maxmin algorithm, random selection, PH landmarks I using $out_{\text{PH}}^{0,1,2}(y)$, PH landmarks II using $out_{\text{PH}}^{1}(y)$, $k--$ landmarks using both cluster centres and outliers as landmarks (kMinusMinus), and $k--$ landmarks using only cluster centres as landmarks (kMinusMinusOutlierFree).}
\end{figure}

\begin{figure}[htp]
\centering
\includegraphics[width=.3\textwidth]{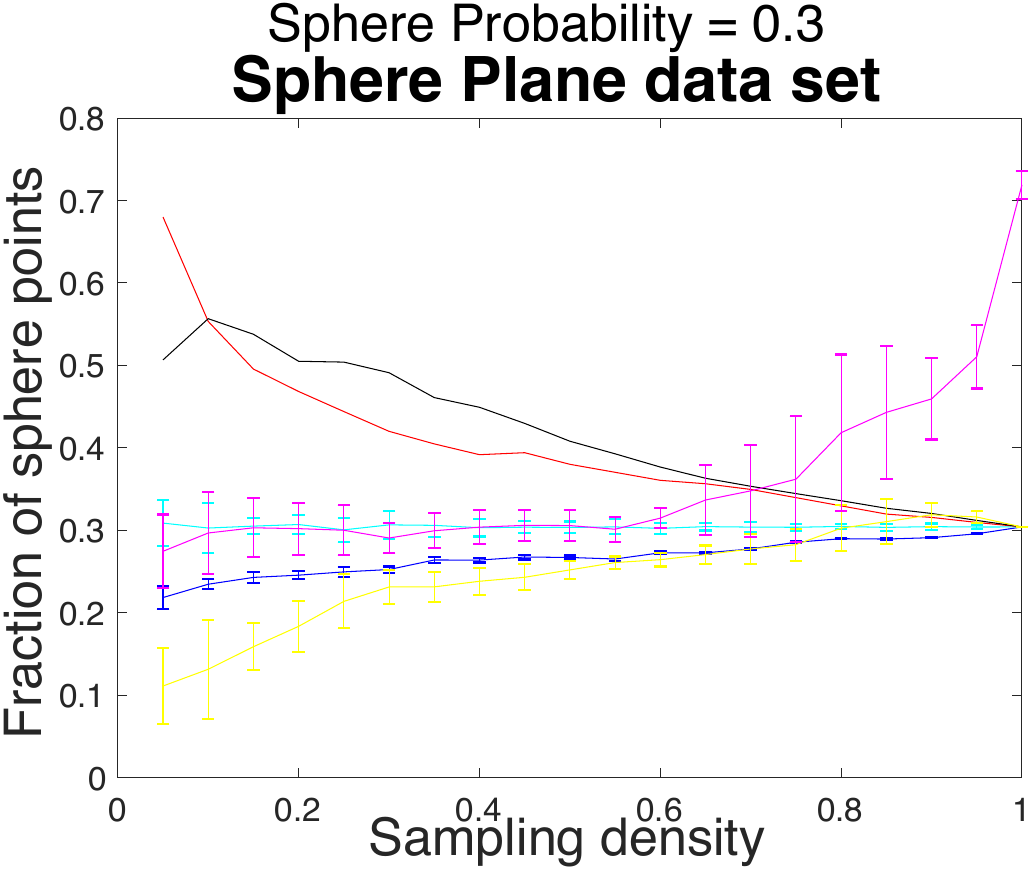}
\includegraphics[width=.3\textwidth]{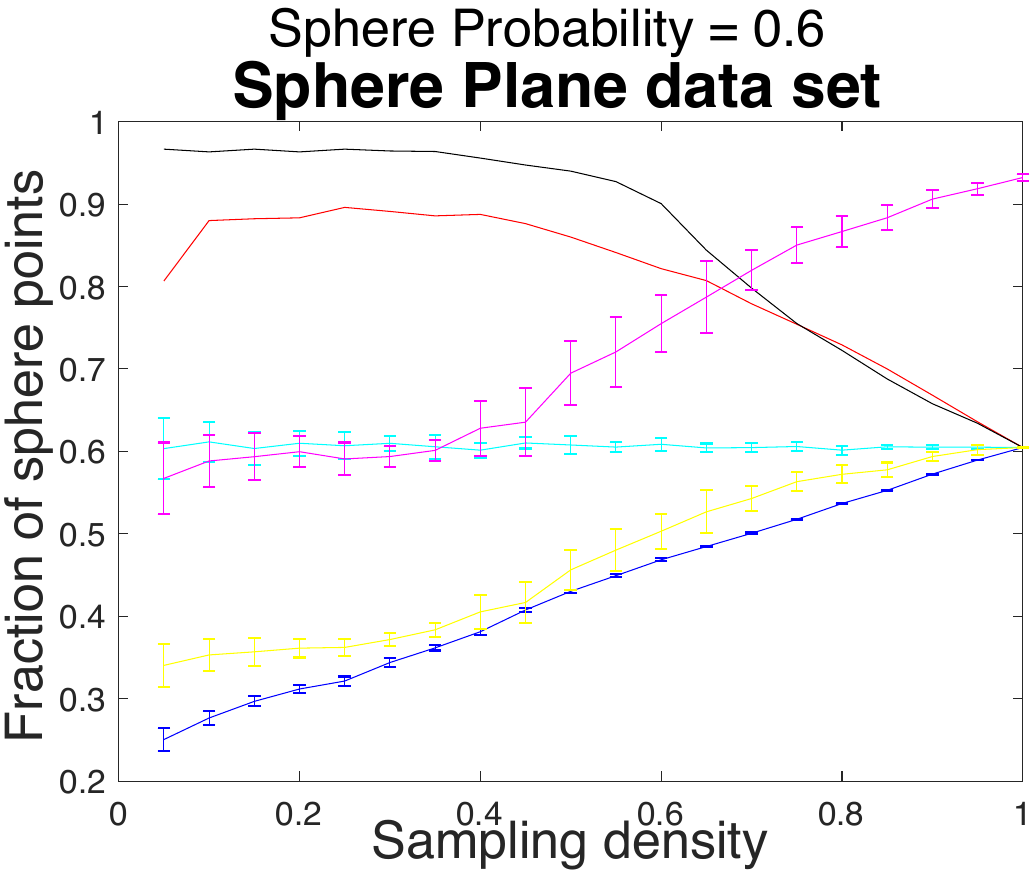}
\includegraphics[width=.3\textwidth]{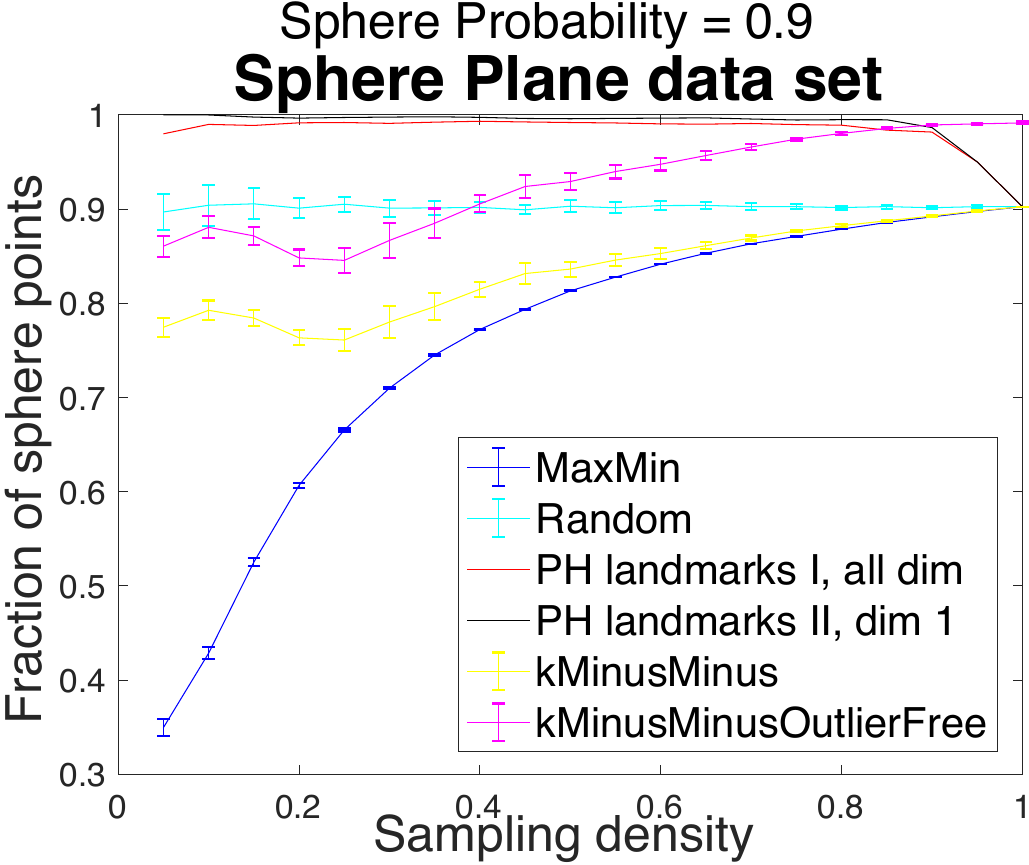}
\caption{\label{ComparisonSpherePlane} Comparison of the fraction of sphere points in selected landmark points for different landmark selection techniques on the sphere-plane dataset for $\delta = 0.2$.We consider landmark selection via the maxmin algorithm, random selection, PH landmarks I using $out_{\text{PH}}^{0,1,2}(y)$, PH landmarks II using $out_{\text{PH}}^{1}(y)$, $k--$ landmarks using both cluster centres and outliers as landmarks (kMinusMinus), and $k--$ landmarks using only cluster centres as landmarks (kMinusMinusOutlierFree).}

\end{figure}

\begin{figure}[htp]
\centering
\includegraphics[width=.3\textwidth]{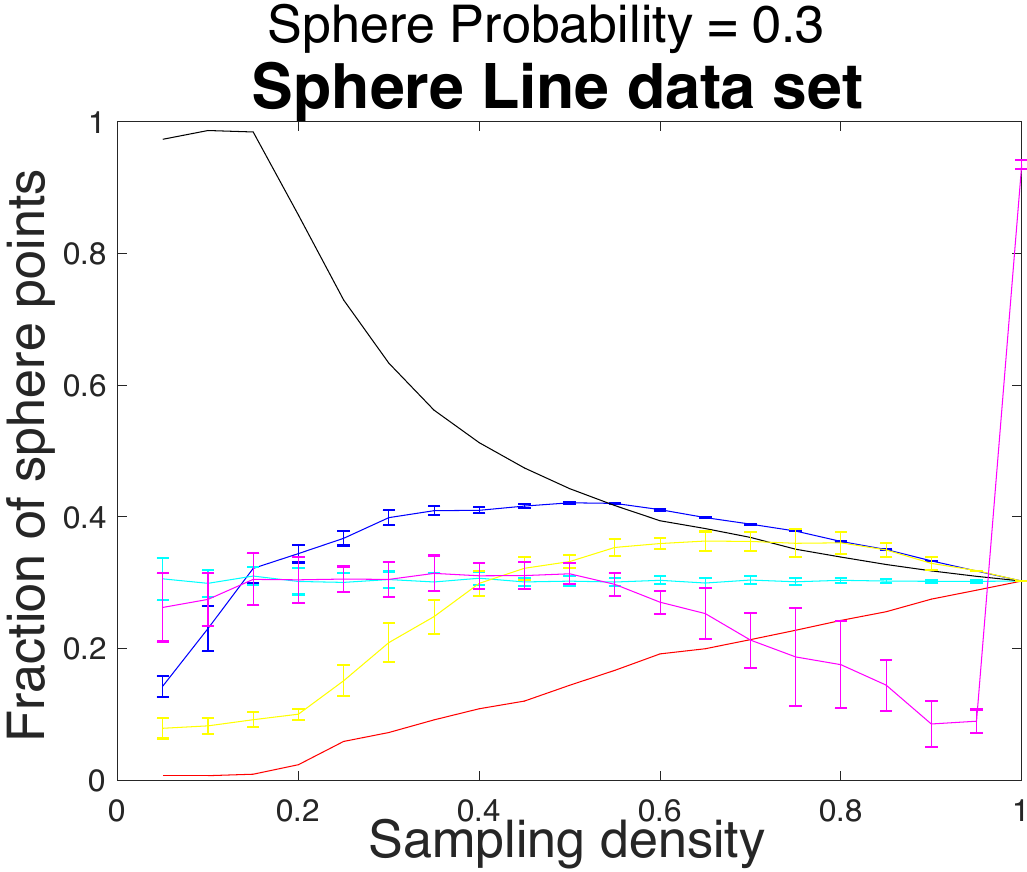}
\includegraphics[width=.3\textwidth]{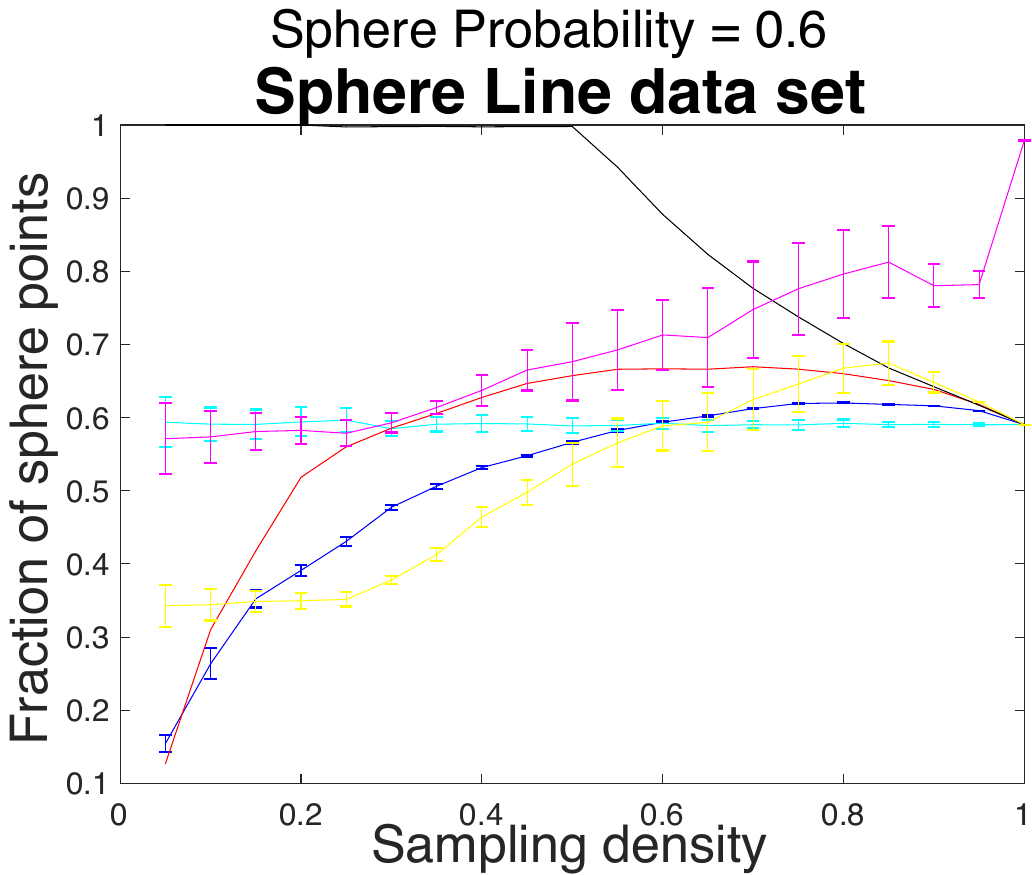}
\includegraphics[width=.3\textwidth]{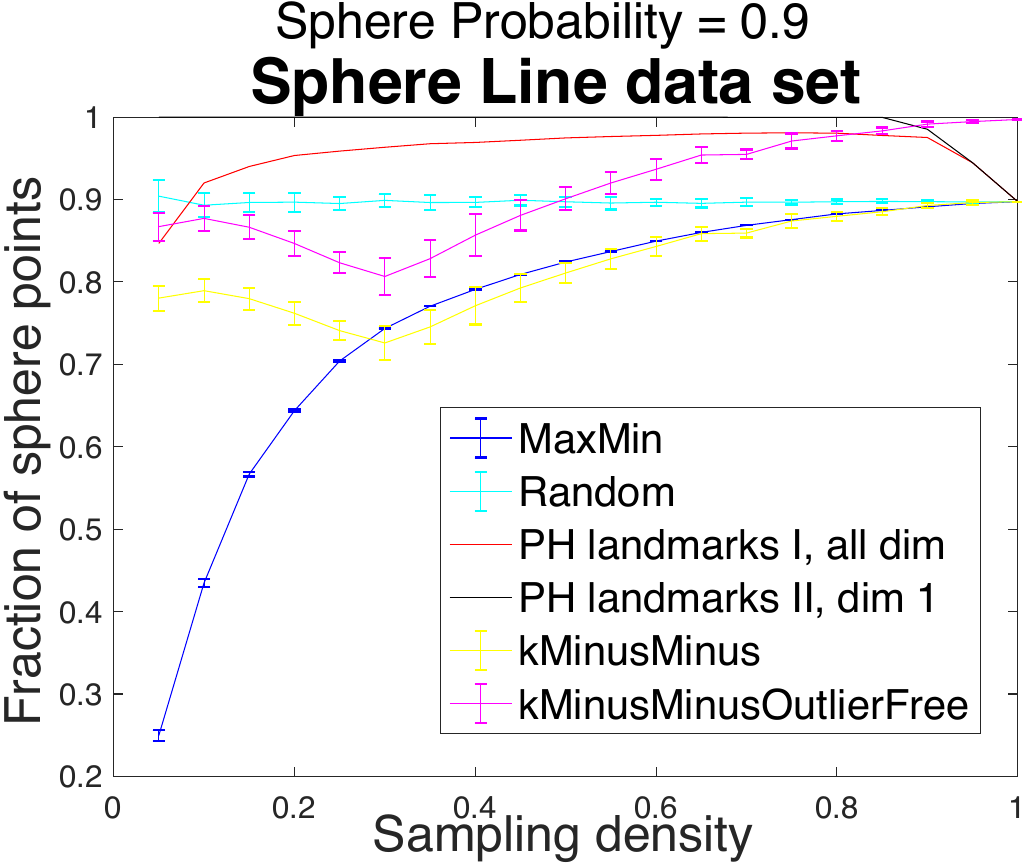}
\caption{\label{ComparisonSphereLine} Comparison of the fraction of sphere points in selected landmark points for different landmark selection techniques on the sphere-line dataset for $\delta = 0.2$.We consider landmark selection via the maxmin algorithm, random selection, PH landmarks I using $out_{\text{PH}}^{0,1,2}(y)$, PH landmarks II using $out_{\text{PH}}^{1}(y)$, $k--$ landmarks using both cluster centres and outliers as landmarks (kMinusMinus), and $k--$ landmarks using only cluster centres as landmarks (kMinusMinusOutlierFree).}
\end{figure}

\begin{figure}[htp]
\centering
\includegraphics[width=.3\textwidth]{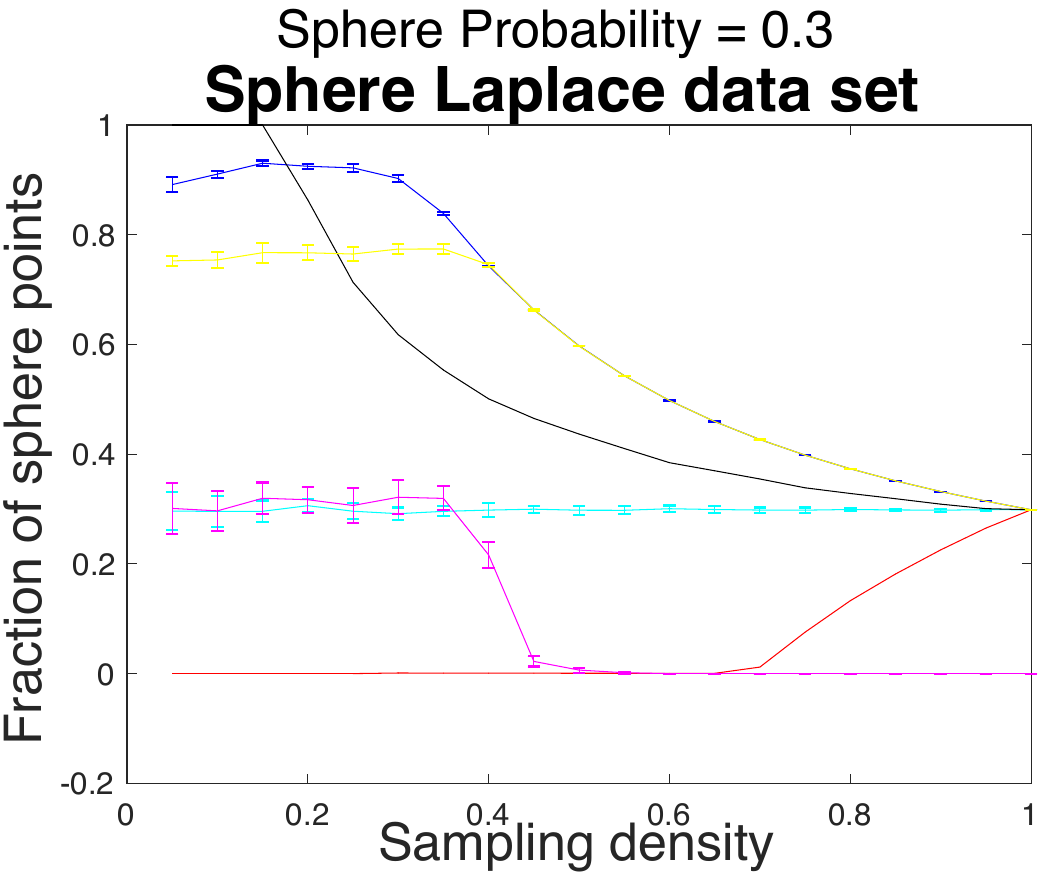}
\includegraphics[width=.3\textwidth]{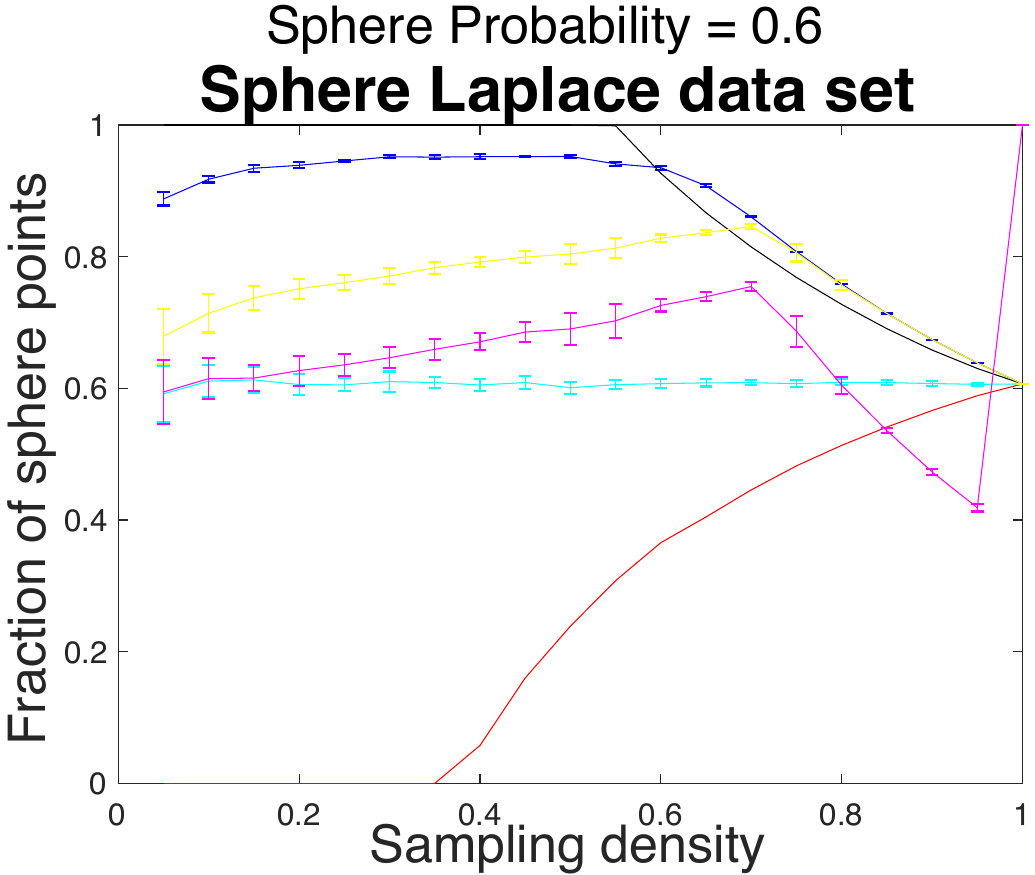}
\includegraphics[width=.3\textwidth]{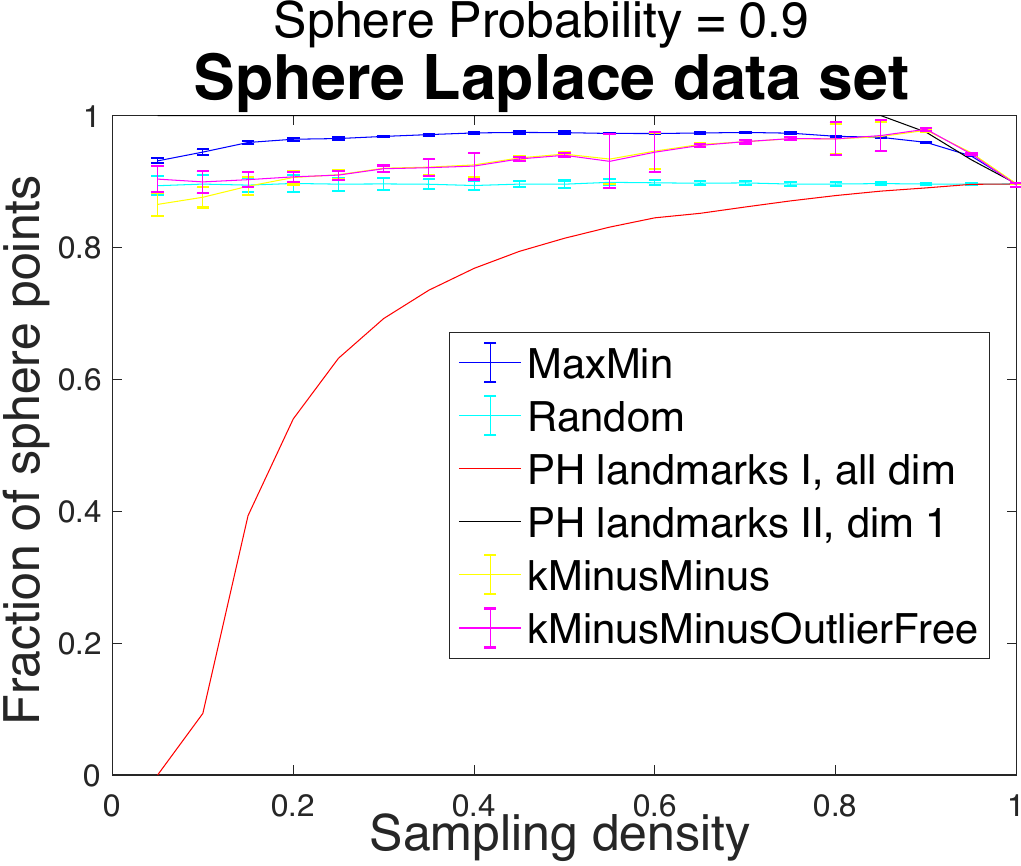}
\caption{\label{ComparisonSphereLaplace} Comparison of the fraction of sphere points in selected landmark points for different landmark selection techniques on the sphere-Laplace dataset for $\delta = 0.2$.We consider landmark selection via the maxmin algorithm, random selection, PH landmarks I using $out_{\text{PH}}^{0,1,2}(y)$, PH landmarks II using $out_{\text{PH}}^{1}(y)$, $k--$ landmarks using both cluster centres and outliers as landmarks (kMinusMinus), and $k--$ landmarks using only cluster centres as landmarks (kMinusMinusOutlierFree).}
\end{figure}

\begin{figure}[htp]
\centering
\includegraphics[width=.3\textwidth]{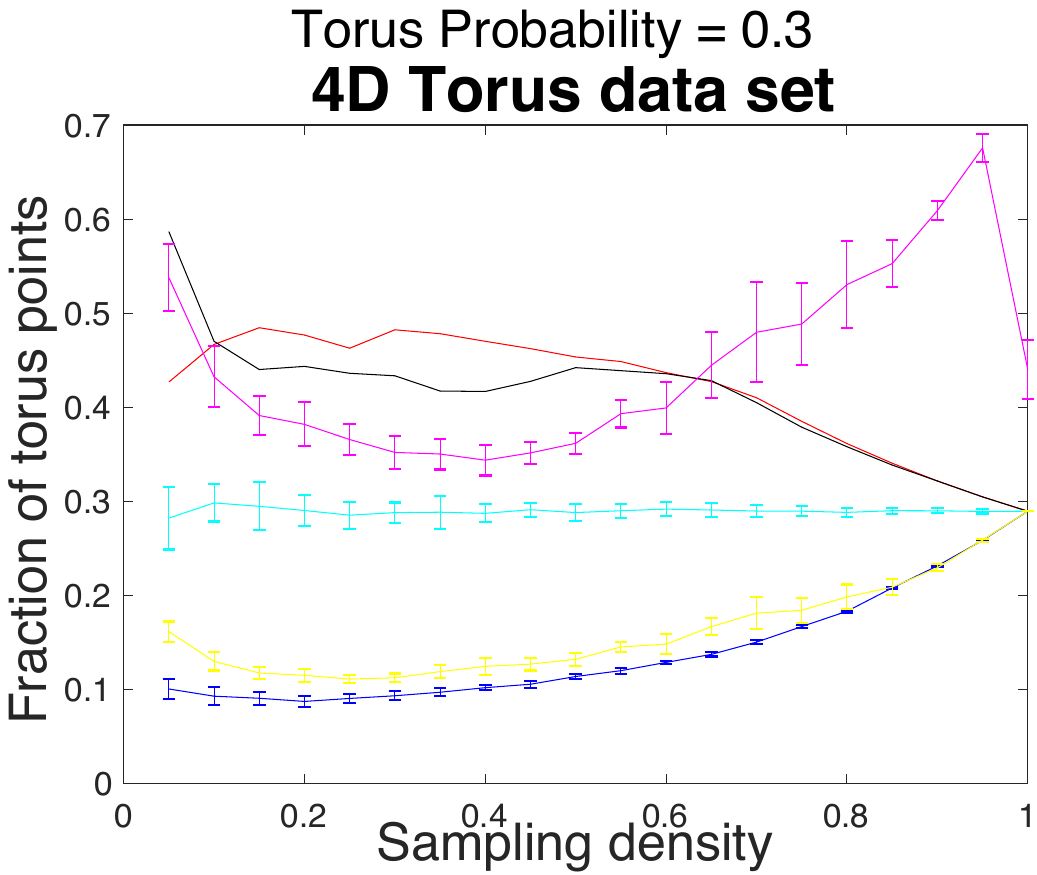}
\includegraphics[width=.3\textwidth]{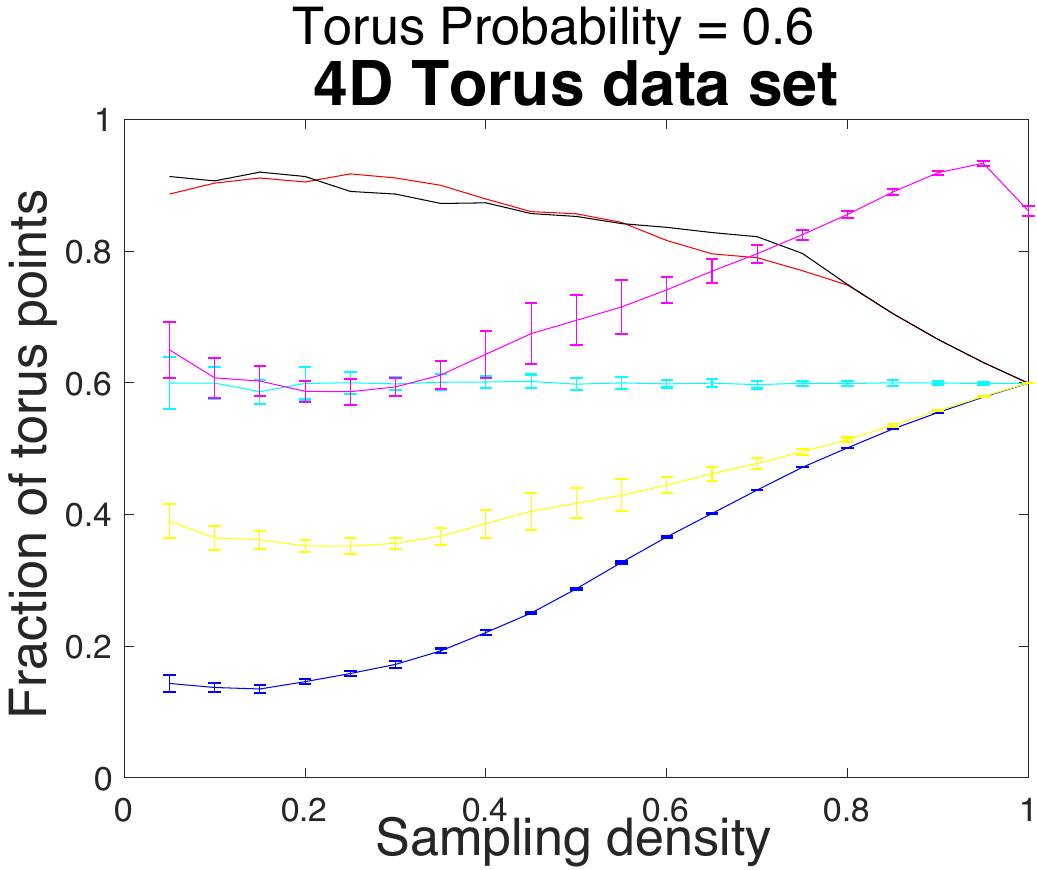}
\includegraphics[width=.3\textwidth]{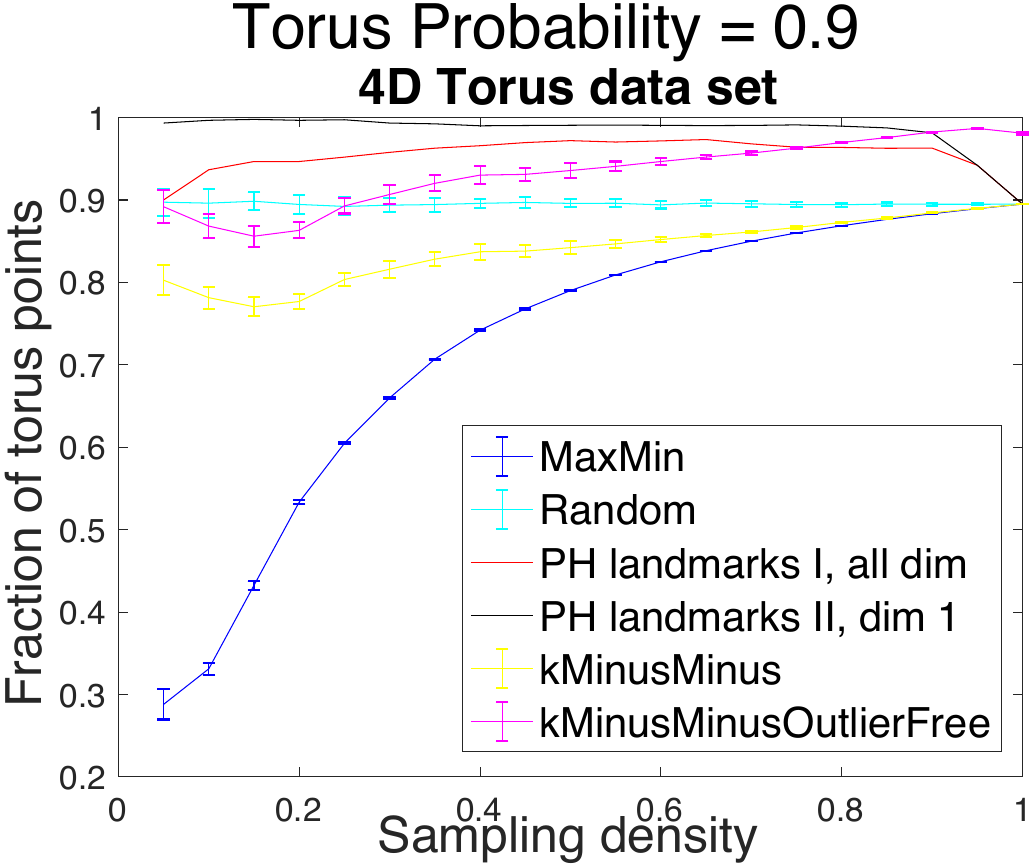}
\caption{\label{ComparisonTorus} Comparison of the fraction of sphere points in selected landmark points for different landmark selection techniques on the Torus dataset for $\delta = 0.5$.We consider landmark selection via the maxmin algorithm, random selection, PH landmarks I using $out_{\text{PH}}^{0,1,2}(y)$, PH landmarks II using $out_{\text{PH}}^{1}(y)$, $k--$ landmarks using both cluster centres and outliers as landmarks (kMinusMinus), and $k--$ landmarks using only cluster centres as landmarks (kMinusMinusOutlierFree).}
\end{figure}

\begin{figure}[htp]
\centering
\includegraphics[width=.3\textwidth]{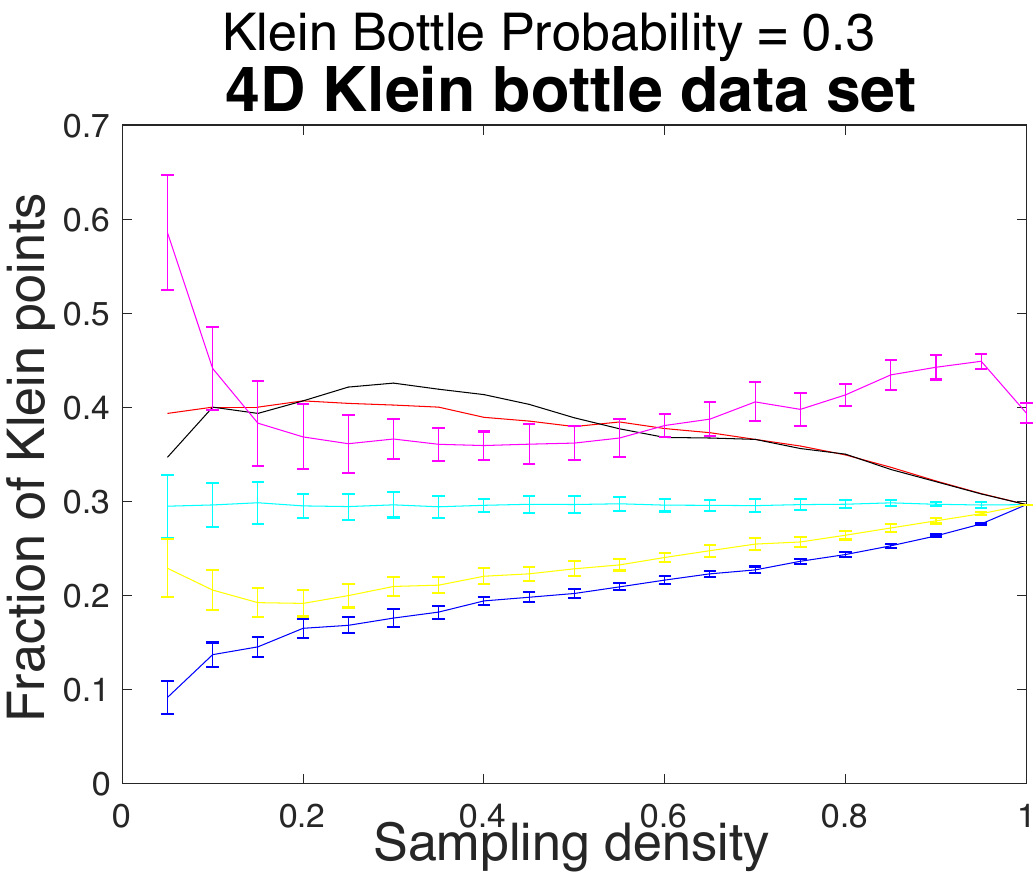}
\includegraphics[width=.3\textwidth]{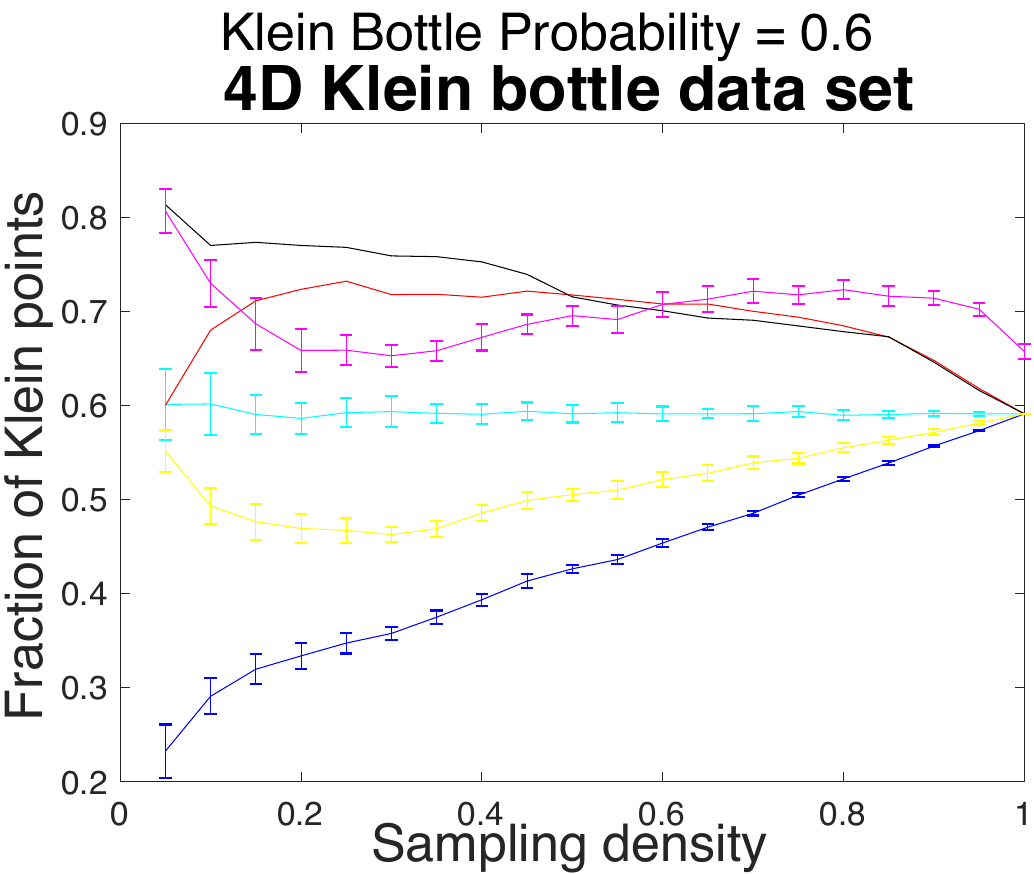}
\includegraphics[width=.3\textwidth]{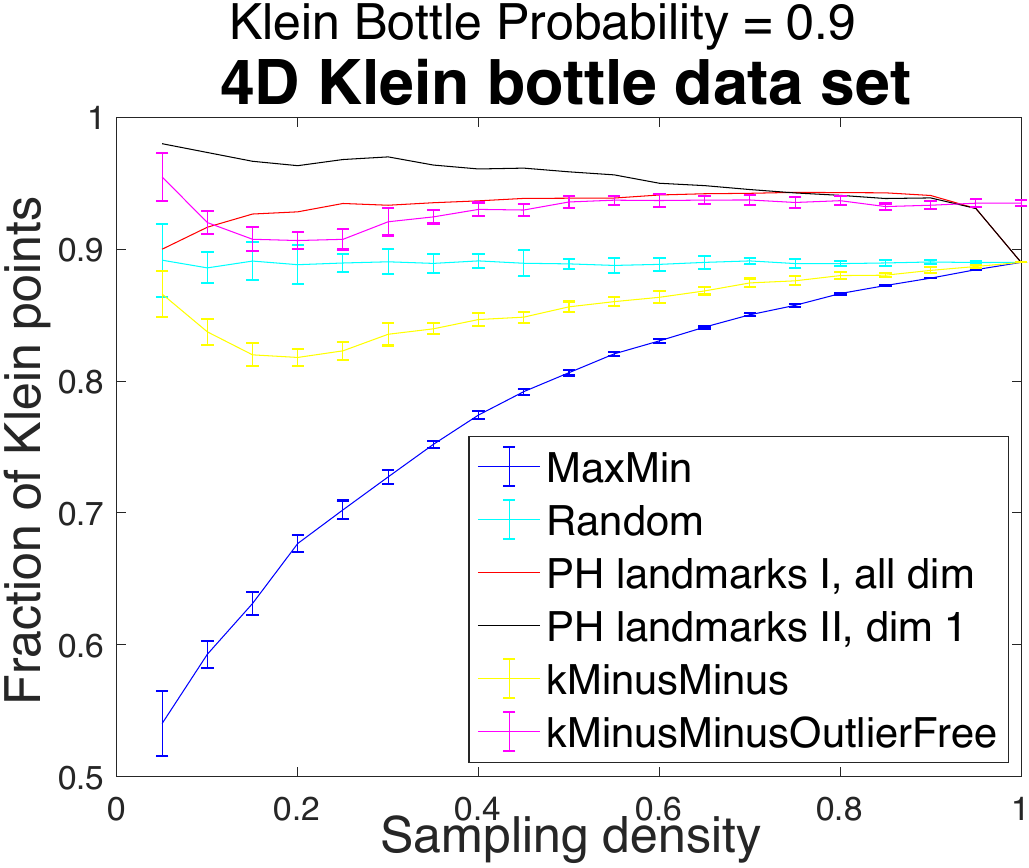}
\caption{\label{ComparisonKlein} Comparison of the fraction of sphere points in selected landmark points for different landmark selection techniques on the Klein bottle dataset for $\delta = 0.6$.We consider landmark selection via the maxmin algorithm, random selection, PH landmarks I using $out_{\text{PH}}^{0,1,2}(y)$, PH landmarks II using $out_{\text{PH}}^{1}(y)$, $k--$ landmarks using both cluster centres and outliers as landmarks (kMinusMinus), and $k--$ landmarks using only cluster centres as landmarks (kMinusMinusOutlierFree).}
\end{figure}

For all our datasets, our aim for the landmarks is to contain a high fraction of signal points, even when sampling only a small fraction of the data as landmarks. In Figures~\ref{ComparisonSphereCube}--\ref{ComparisonKlein} we show plots of the fraction of signal points in the various landmark sets at different sampling densities. Since in the PH landmark selection we allow super outliers as landmarks once all other points are taken, we expect the fraction of signal landmarks for sampling density 1 to represent the probability of signal points in the dataset, except in the variant of the $k--$ landmarks where we include only the cluster centres as landmarks (referred to as `kMinusMinusOutlierFree' in the plots). Note that for this variant of $k--$ landmarks for datasets with $p\leq 0.5$, it is possible to obtain a signal fraction of $0$ even for sampling density 1 if the algorithm selects all $k = p m$ cluster centres to be located among noise points. For the maxmin, random and $k--$ landmarks we show the average fraction of signal points and its standard deviation across 20 realisations of the selection algorithms.

As expected, we observe that the  maxmin algorithm tends to select noise points as landmarks for all datasets with only one exception (see Figures~\ref{ComparisonSphereCube}--\ref{ComparisonKlein}): for the sphere-Laplace data, the  maxmin algorithm performs well (see Fig.~\ref{ComparisonSphereLaplace}), as the noise is located in a cluster far away from the signal and hence maximising the distance between landmarks results in many points being selected from the sphere. The fraction of signal landmarks for random selection also behaves as we expect for all datasets in Figures~\ref{ComparisonSphereCube}--\ref{ComparisonKlein}: the selected landmarks are representative for the whole dataset with an almost constant fraction of signal points over all sampling densities that corresponds to the fraction of signal points in the dataset. For the $k--$ landmarks, we can see a clear improvement in the signal fraction in most datasets in Figures~\ref{ComparisonSphereCube}--\ref{ComparisonKlein} when considering only cluster centres as landmarks -- the inclusion of outliers gives the $k--$ landmarks similarly bad properties as the  maxmin landmarks. The $k--$ landmarks that do not include outliers tend to perform well for high sampling densities, where the number outlier points corresponds roughly to the number of noise points in the dataset. For low sampling densities however, the method only outperforms random selection for most of the sphere-cube and Klein bottle datasets (see Figures~\ref{ComparisonSphereCube} and \ref{ComparisonKlein}). For the sphere-plane, sphere-line and the sphere Laplace-line datasets the reason for this lies in the nature of the noise, which leads to the selection of cluster centres in the noise data. We also notice large standard deviations from the average fraction of signal points in the $k--$ landmark over 20 realisations.

With exception of the sphere-line and sphere-Laplace datasets (see Figures~\ref{ComparisonSphereLine} and \ref{ComparisonSphereLaplace}), the PH landmark selection techniques I and II both outperform the standard methods as well as the $k--$ landmarks clearly for most cases, especially for low sampling densities. Interestingly, the $k--$ landmarks perform very well on the Klein bottle dataset (see Fig.~\ref{ComparisonKlein}), beating both PH landmarks for very small sampling densities.
For the sphere-line and sphere-Laplace datasets (see Figures~\ref{ComparisonSphereLine} and \ref{ComparisonSphereLaplace}), PH landmarks II restricted to dimension 1 outperform all other methods for low sampling densities while PH landmarks I across all dimensions perform worse than all other methods in most cases. The noise in these datasets is located on lines in dense regions of the dataset where the local PH does not find any topological features in dimension 1, but many features with low persistence in dimension 0.
In general, the two versions of PH landmarks start coinciding as soon as super outliers are added to the dataset which we can observe in the plots as a rapid drop in the fraction of signal points. We add the super outliers to the landmarks in random order (once all other points are already selected as landmarks) and hence both PH landmark methods differ only slightly in the development of their signal fractions after the addition of super outliers. 
We note there are cases in which the PH landmarks thrive because the points that are not super outliers are predominantly signal points. 
This is not the case for the sphere-line and the sphere-Laplace datasets with $p = 0.6$ (see Figures~\ref{ComparisonSphereLine} and \ref{ComparisonSphereLaplace}) for which the dimension~1 PH landmarks II perform very well. Both of these datasets have less than 4 super outliers. Interestingly, there seems to be a trend for the dimension 1 PH landmarks II to outperform PH landmarks I on datasets with high signal content, i.e. for $p > 0.6$ across all datasets. For lower signal content the PH landmarks I perform strongly. For the Klein bottle dataset, dimension~1 PH landmarks II outperform PH landmarks I in almost all cases.

Overall, the results underline that both PH landmarks I and II represent the PH of the dataset well and are robust to outliers, in particular for low sampling densities. They outperform standard methods in a large majority of cases and, moreover, via the PH outlierness score they give us a notion of how much a point can be considered an outlier for the respective variant of the method. $k--$ landmarks perform better than random selection in most cases and perform very well, in particular, for high sampling densities. Given the much higher computational cost and the large fluctuations in signal fraction between different realisations of the method, using random selection instead of $k--$ landmarks could however present a more practical approach.

\subsection{Comparisons between persistent homology landmark selection methods and dense core subsets}\label{Subsubsec:DenseCoreComparison}

\begin{figure}[htp]
\centering
\includegraphics[width=.49\textwidth]{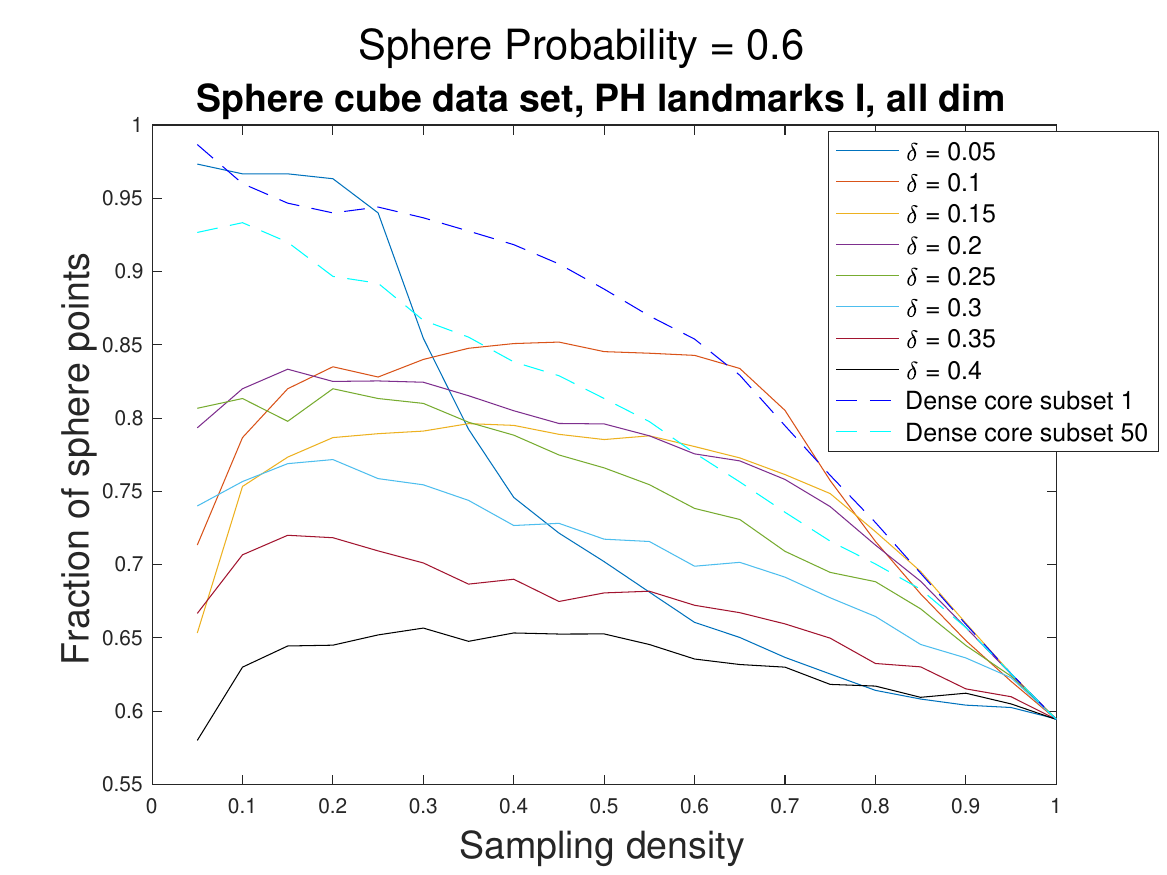}
\includegraphics[width=.49\textwidth]{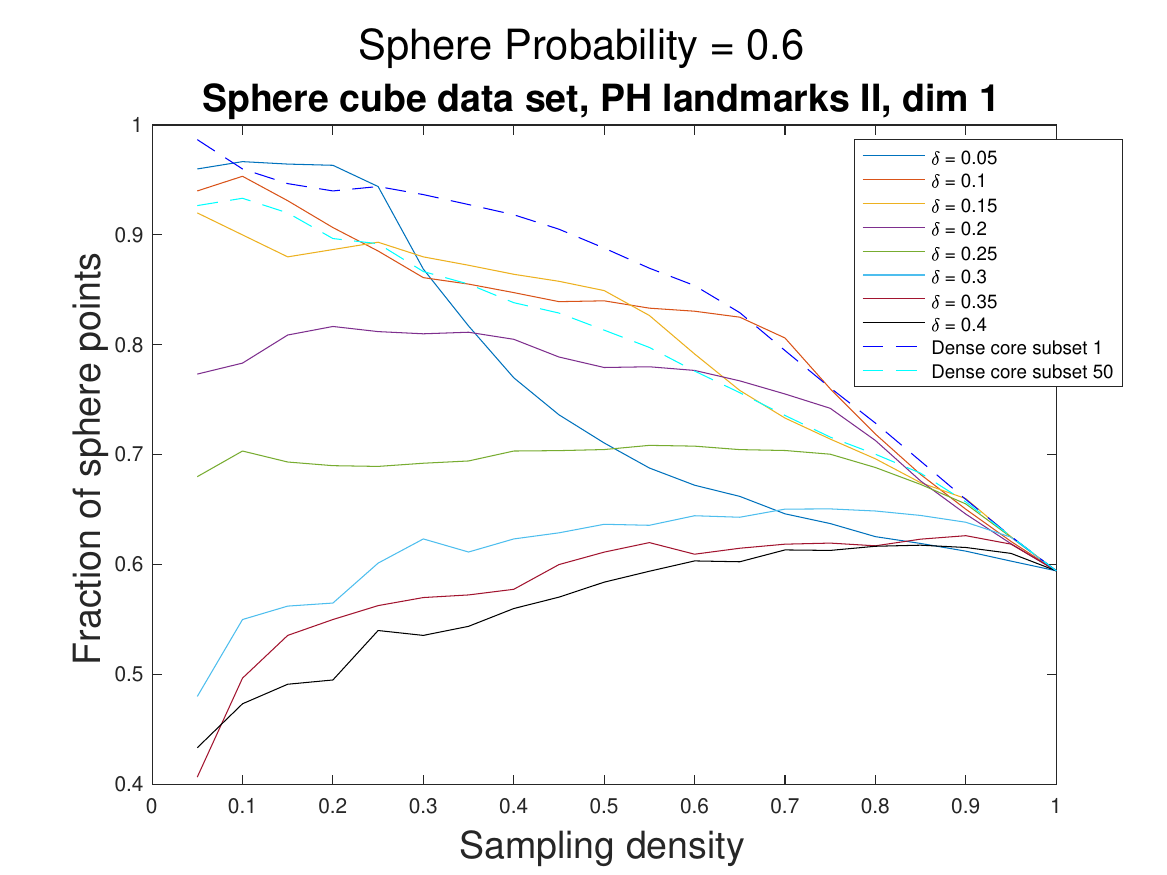}
\caption[Comparison of the fraction of sphere points in selected PH landmark points for different values of $\delta$ and dense core subsets on the sphere-cube dataset, $p = 0.6$.]{\label{ComparisonSphereCubeCore} Comparison of the fraction of sphere points in selected PH landmark points for different values of $\delta$ and dense core subsets on the sphere-cube dataset, $p = 0.6$. We show PH landmarks I using $out_{\text{PH}}^{0,1,2}(y)$ (left) and PH landmarks II using $out_{\text{PH}}^{1}(y)$ (right).}
\end{figure}

\begin{figure}[htp]
\centering
\includegraphics[width=.49\textwidth]{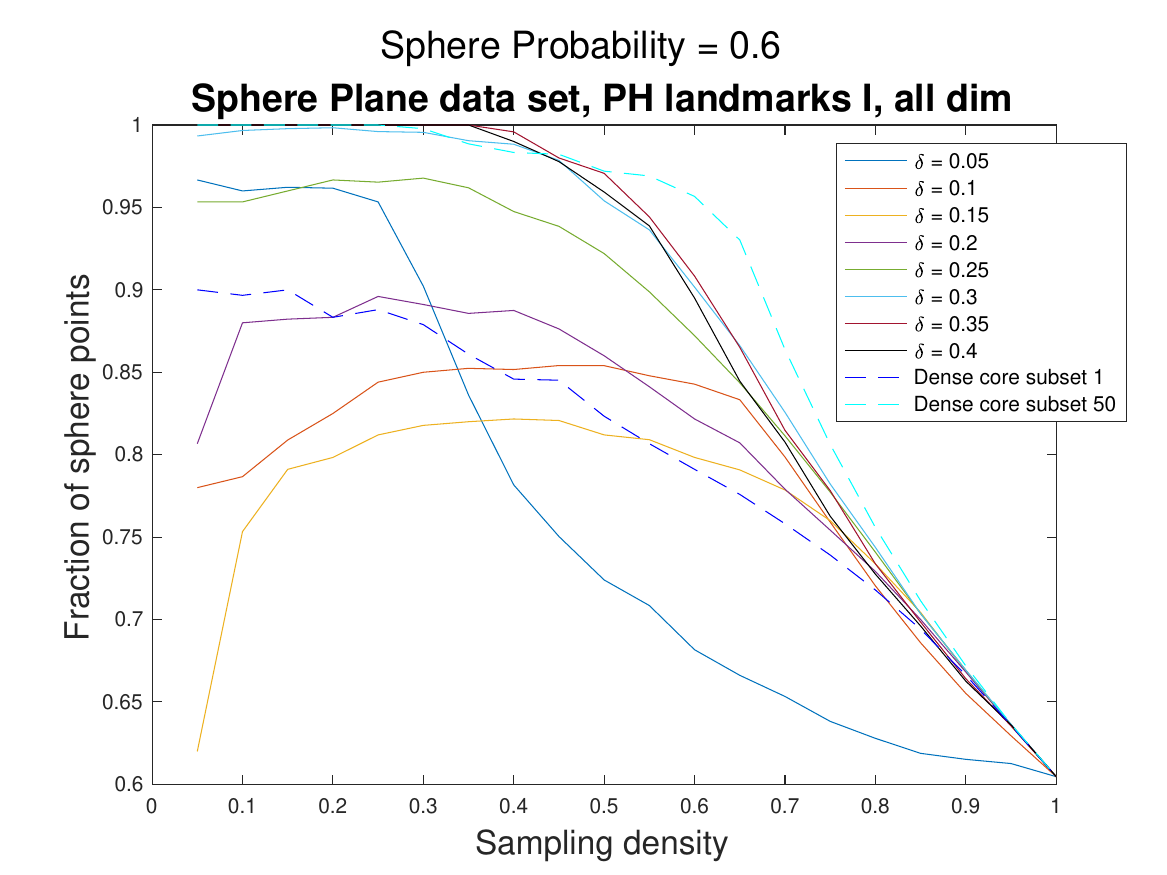}
\includegraphics[width=.49\textwidth]{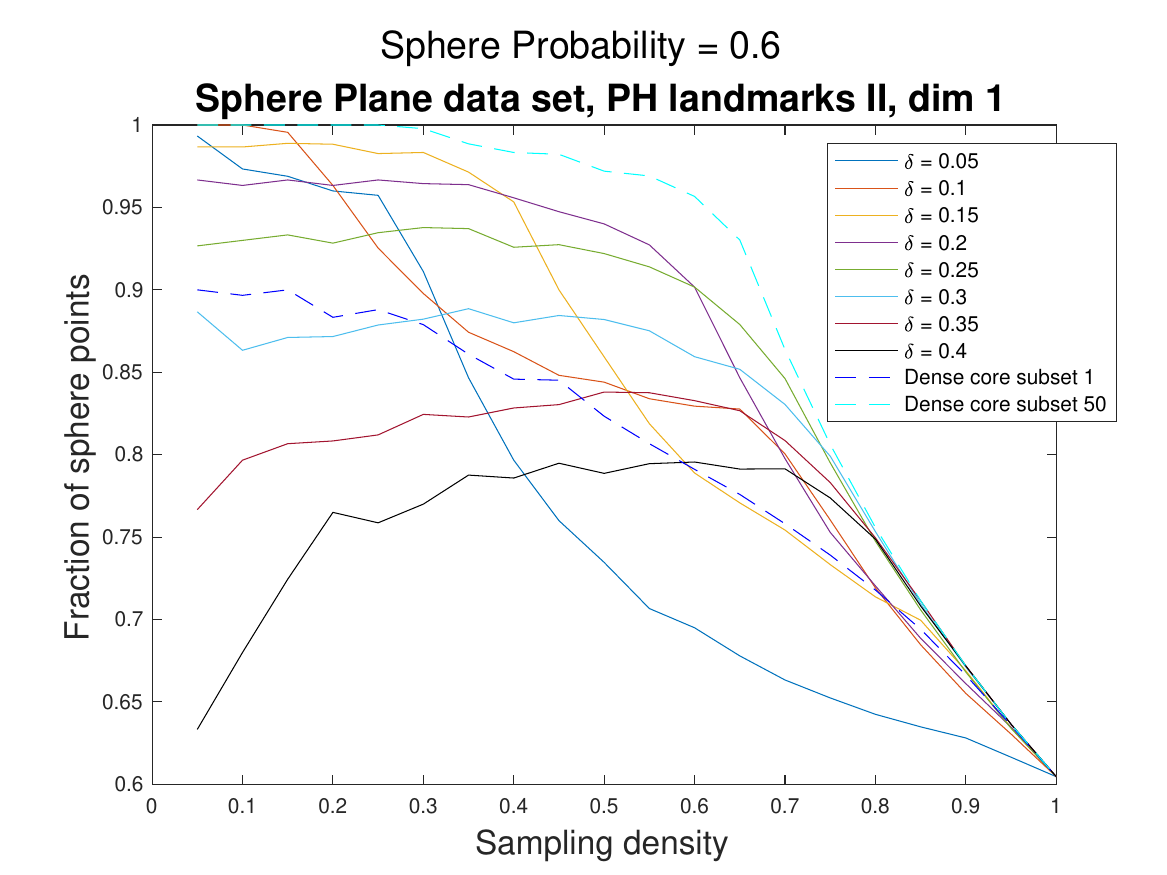}
\caption[Comparison of the fraction of sphere points in selected PH landmark points for different values of $\delta$ and dense core subsets on the sphere-plane dataset, $p = 0.6$.]{\label{ComparisonSpherePlaneCore} Comparison of the fraction of sphere points in selected PH landmark points for different values of $\delta$ and dense core subsets on the sphere-plane dataset, $p = 0.6$. We show PH landmarks I using $out_{\text{PH}}^{0,1,2}(y)$ (left) and PH landmarks II using $out_{\text{PH}}^{1}(y)$ (right).}
\end{figure}

\begin{figure}[htp]
\centering
\includegraphics[width=.49\textwidth]{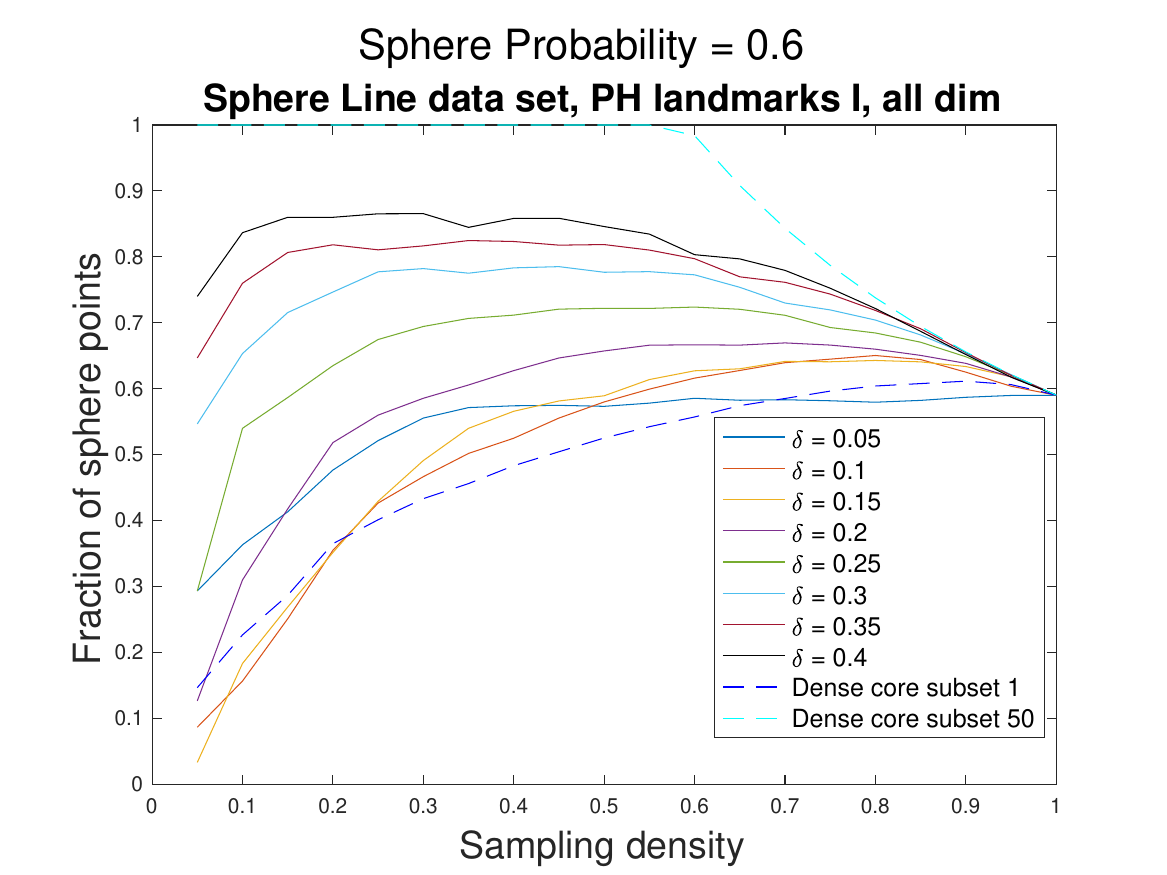}
\includegraphics[width=.49\textwidth]{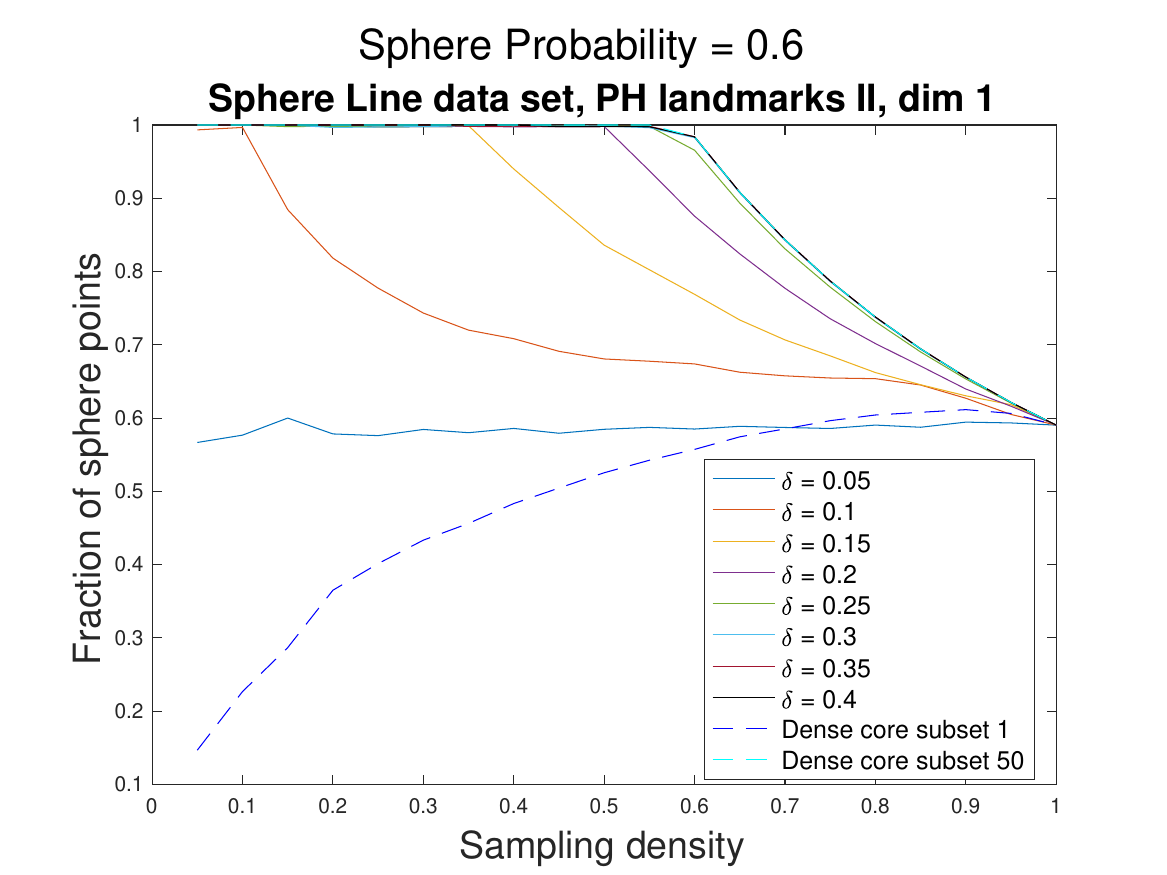}
\caption[Comparison of the fraction of sphere points in selected PH landmark points for different values of $\delta$ and dense core subsets on the sphere-line dataset, $p = 0.6$.]{\label{ComparisonSphereLineCore} Comparison of the fraction of sphere points in selected PH landmark points for different values of $\delta$ and dense core subsets on the sphere-line dataset, $p = 0.6$. We show PH landmarks I using $out_{\text{PH}}^{0,1,2}(y)$ (left) and PH landmarks II using $out_{\text{PH}}^{1}(y)$ (right).}
\end{figure}

\begin{figure}[htp]
\centering
\includegraphics[width=.49\textwidth]{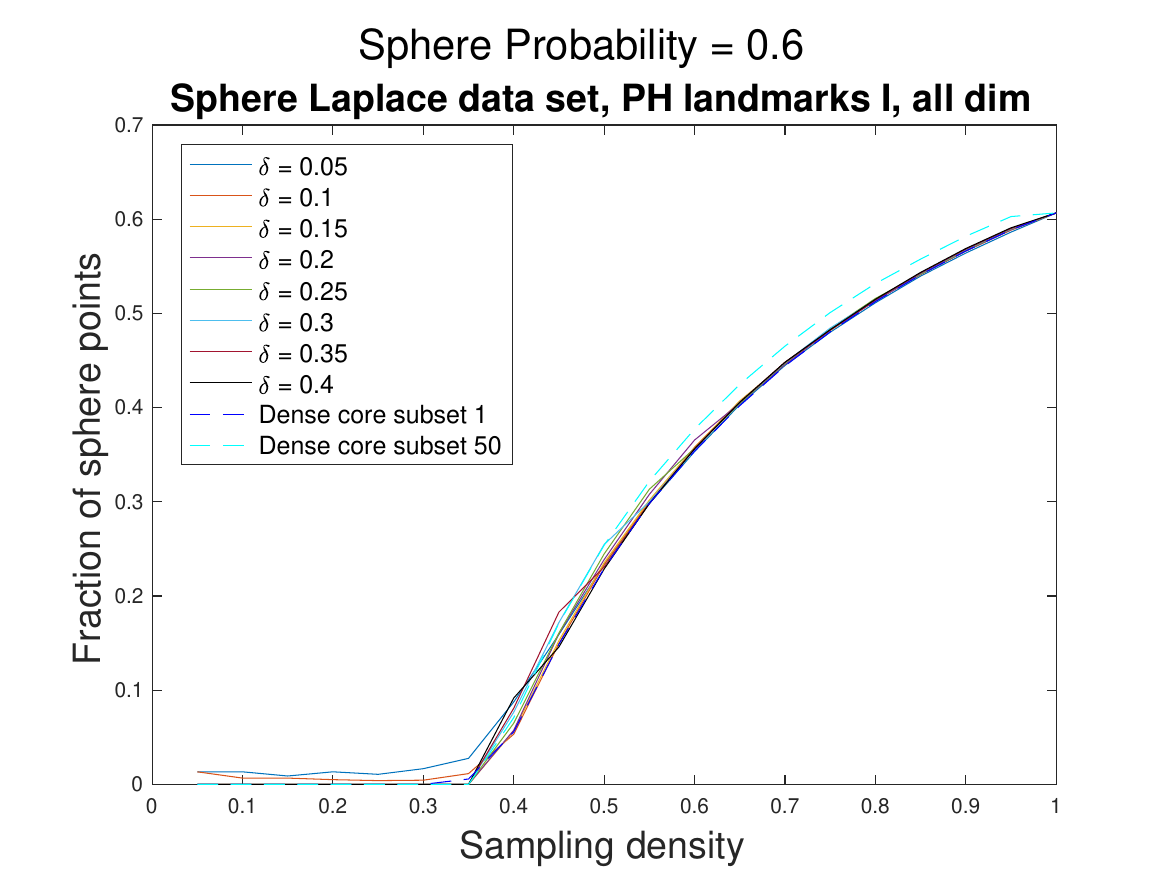}
\includegraphics[width=.49\textwidth]{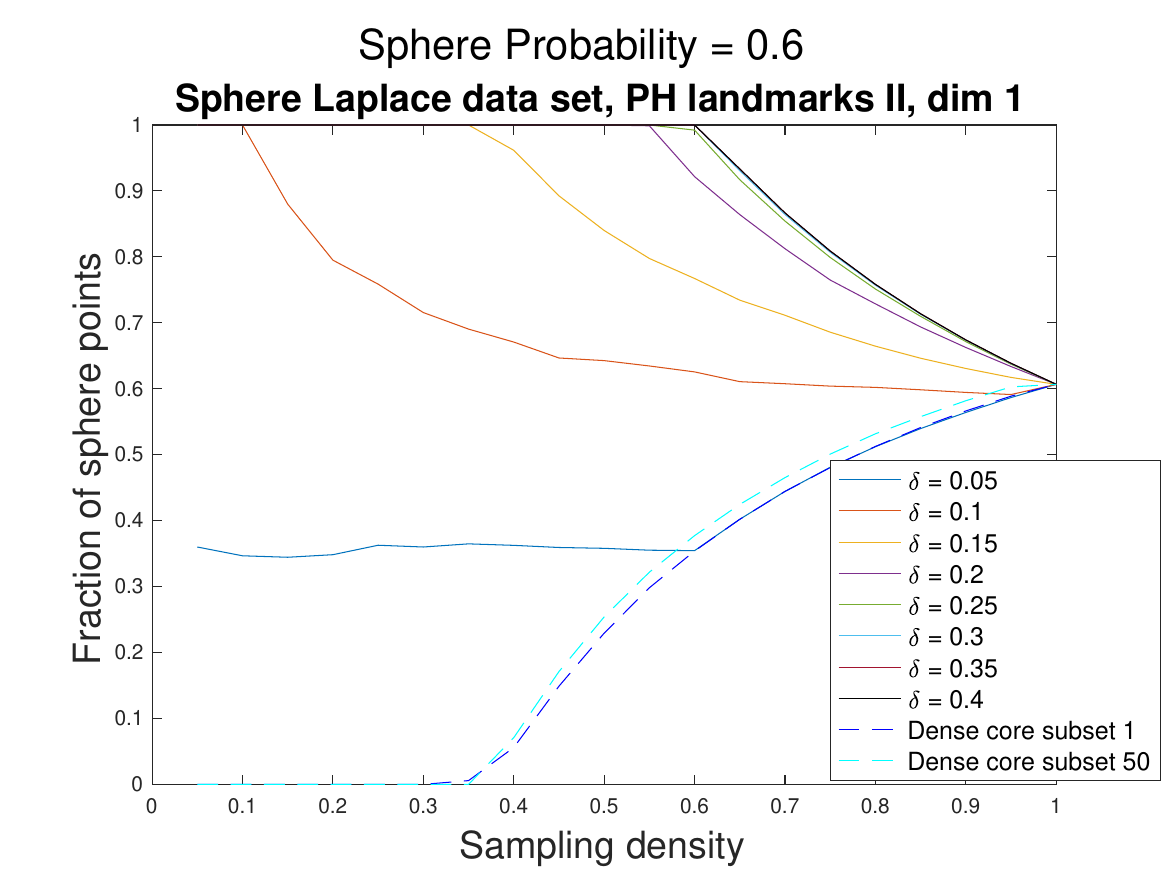}
\caption[Comparison of the fraction of sphere points in selected PH landmark points for different values of $\delta$ and dense core subsets on the sphere-Laplace dataset, $p = 0.6$.]{\label{ComparisonSphereLaplaceCore} Comparison of the fraction of sphere points in selected PH landmark points for different values of $\delta$ and dense core subsets on the sphere-Laplace dataset, $p = 0.6$. We show PH landmarks I using $out_{\text{PH}}^{0,1,2}(y)$ (left) and PH landmarks II using $out_{\text{PH}}^{1}(y)$ (right).}
\end{figure}

\begin{figure}[htp]
\centering
\includegraphics[width=.49\textwidth]{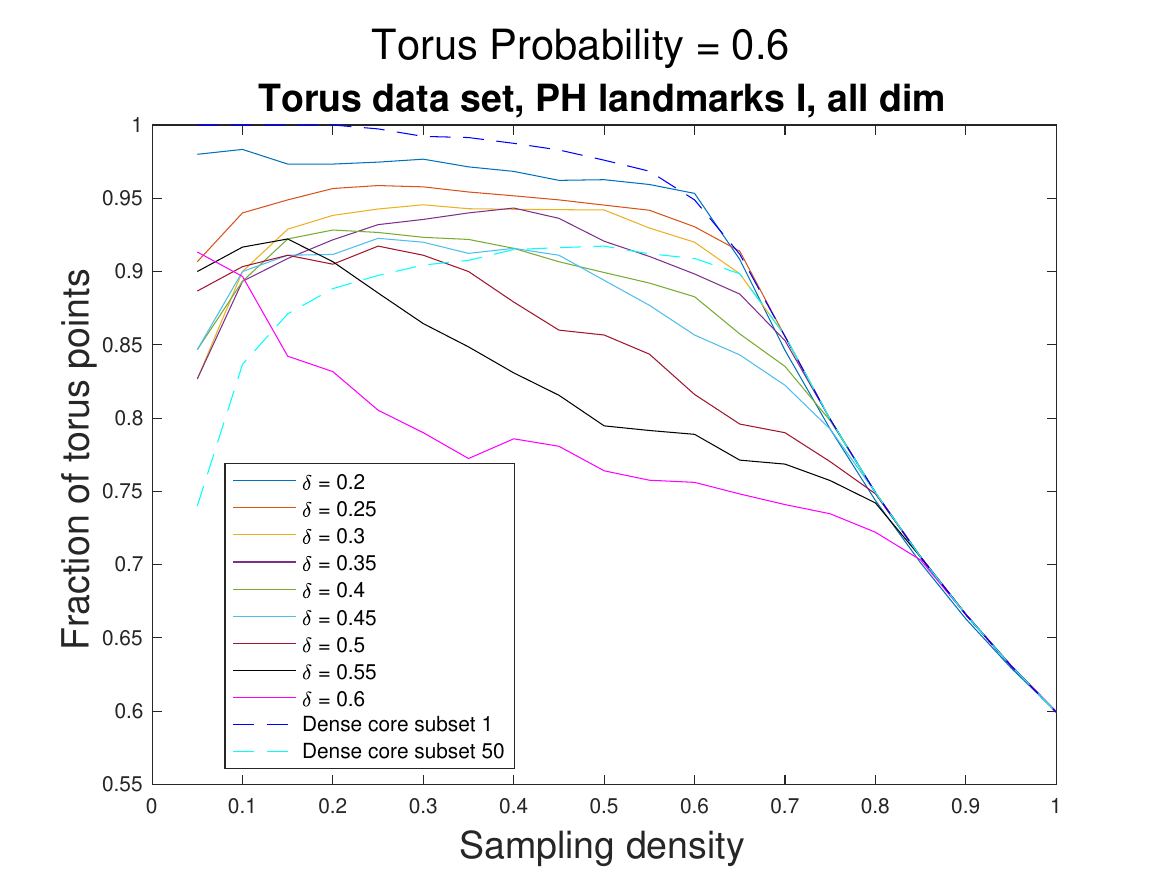}
\includegraphics[width=.49\textwidth]{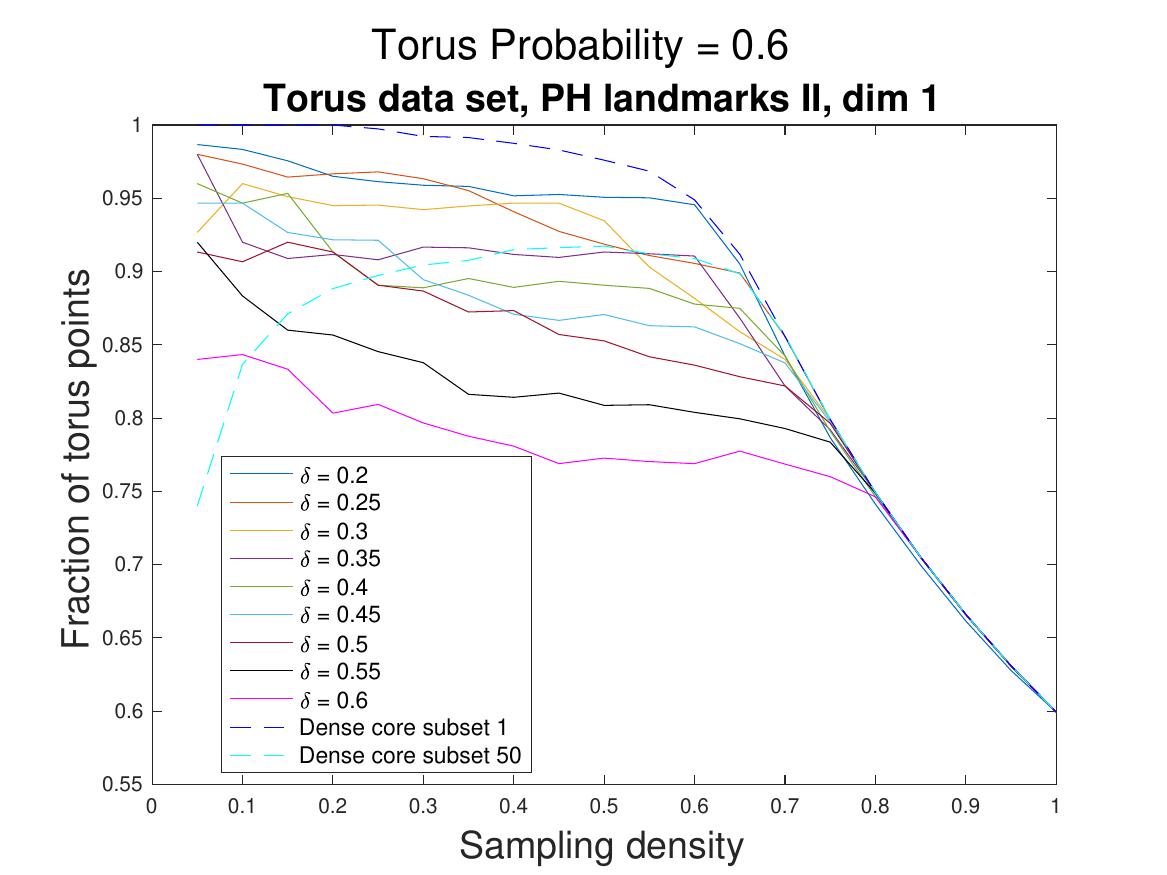}
\caption[Comparison of the fraction of sphere points in selected PH landmark points for different values of $\delta$ and dense core subsets on the torus dataset, $p = 0.6$.]{\label{ComparisonTorusCore} Comparison of the fraction of sphere points in selected PH landmark points for different values of $\delta$ and dense core subsets on the torus dataset, $p = 0.6$. We show PH landmarks I using $out_{\text{PH}}^{0,1,2}(y)$ (left) and PH landmarks II using $out_{\text{PH}}^{1}(y)$ (right).}
\end{figure}

\begin{figure}[htp]
\centering
\includegraphics[width=.49\textwidth]{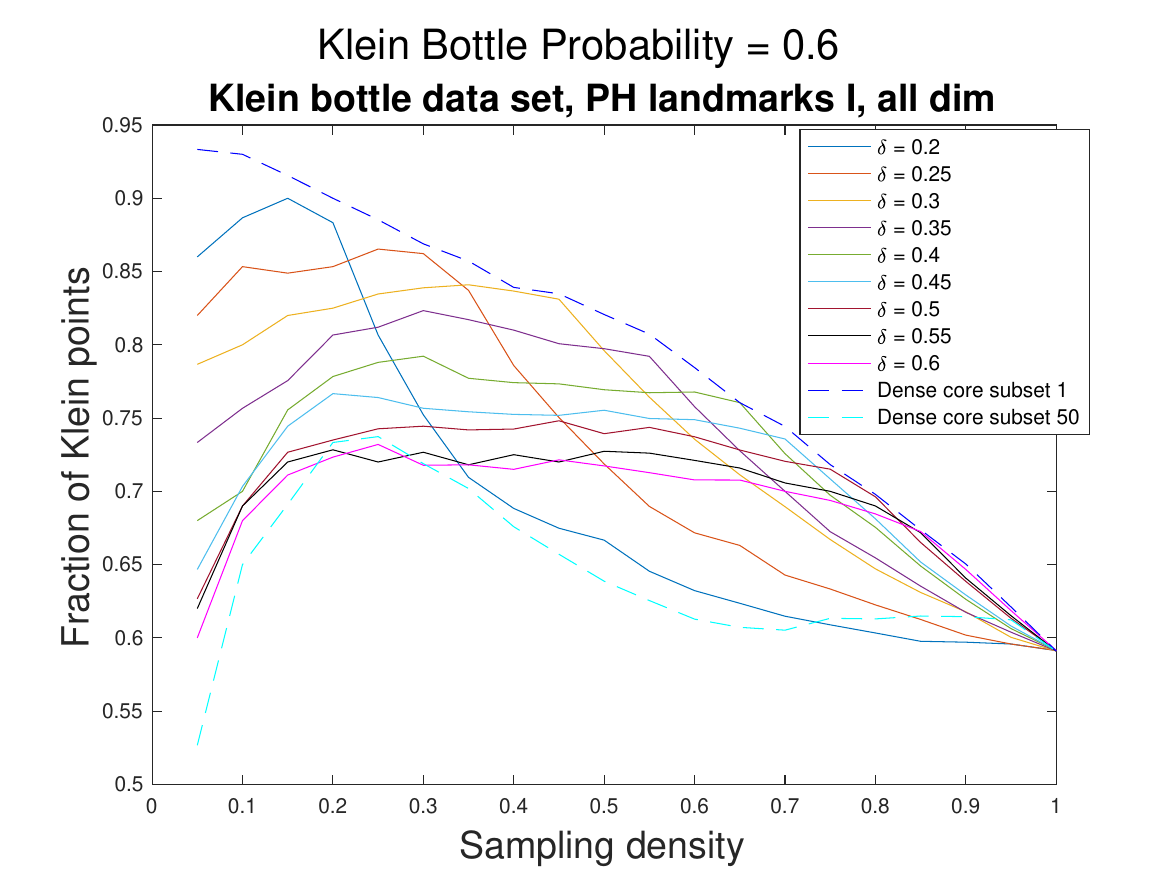}
\includegraphics[width=.49\textwidth]{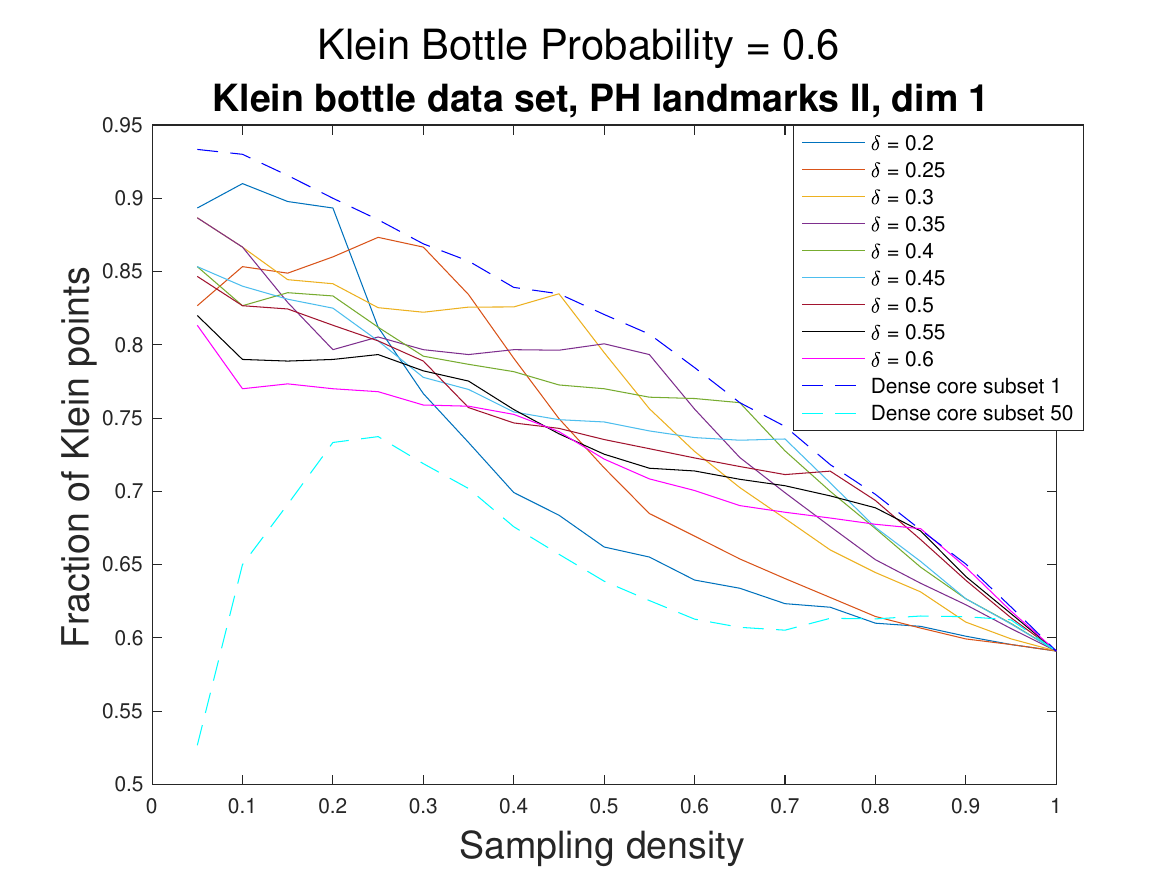}
\caption[Comparison of the fraction of sphere points in selected PH landmark points for different values of $\delta$ and dense core subsets on the Klein bottle dataset, $p = 0.6$.]{\label{ComparisonKleinCore} Comparison of the fraction of sphere points in selected PH landmark points for different values of $\delta$ and dense core subsets on the Klein bottle dataset, $p = 0.6$. We show PH landmarks I using $out_{\text{PH}}^{0,1,2}(y)$ (left) and PH landmarks II using $out_{\text{PH}}^{1}(y)$ (right).}
\end{figure}

We now provide a more in detail study of the two PH landmark selection techniques, in particular concerning the influence of the $\delta$ parameter. We also compare the techniques to two dense core subsets, using for $K = 1$ and $K = 50$, which we consider as landmark sets. We present our results in Figures~\ref{ComparisonSphereCubeCore} -- \ref{ComparisonKleinCore}. We present only the plots for datasets with a signal probability $p = 0.6$. 

For the dense core subsets, we find that for all datasets except the sphere-plane, sphere-line, and sphere-Laplace datasets (see Figures \ref{ComparisonSpherePlaneCore} -- \ref{ComparisonSphereLaplaceCore}) the local density measure $K = 1$ outperforms the more global density measure $K = 50$. Indeed, in these cases, the dense core subset with $K=1$ captures a larger fraction of signal points than most of our PH landmarks. For the sphere-cube dataset (see Fig.~\ref{ComparisonSphereCubeCore}), we seem to outperform the dense core subset with $K=1$ for a small range of low sampling densities for $\delta = 0.05$, which is a $\delta$ value where most of the data points are classified as super outliers. Our definition of super outliers as points with less than two neighbours in their $\delta$-neighbourhoods, seems to imply that, for this dataset, the distance to the second closest neighbour is more relevant for low sampling densities than the distance to the closest neighbour. Interestingly, for the datasets where the local dense core subset with $K=1$ performs well (see Figures~\ref{ComparisonSphereCubeCore}, \ref{ComparisonTorusCore}, and \ref{ComparisonKleinCore} ), we also see that smaller values of $\delta$ give us better results, both for PH landmarks I across all dimensions and PH landmarks II restricted to dimension 1. For the sphere-line dataset (see Fig.~\ref{ComparisonSphereLineCore}), where we have a better performance for the dense core subset with $K=50$, larger $\delta$-neighbourhoods are more advantageous for both PH landmark versions. In this case, we perform as well as the dense core subset with $K=50$ for $\delta = 0.3, 0.35, 0.4$ for dimension1 PH landmarks II. In the sphere-Laplace dataset (see Fig.~\ref{ComparisonSphereLaplaceCore}) PH landmarks II clearly outperform both dense core subsets. Here again, we find the trend that larger $\delta$ values are advantageous for PH landmarks II. The fact that dense core subsets perform well on most of our datasets is determined by our signal points lying in denser regions of the data than the noise points. It is only in the sphere-Laplace dataset (see Fig.~\ref{ComparisonSphereLaplaceCore}), where the noise does not obey this characteristic and we observe that the local PH information in this case is richer than the distance to the $K$-th neighbour.

\begin{figure}[htp]
\centering
\includegraphics[width=.49\textwidth]{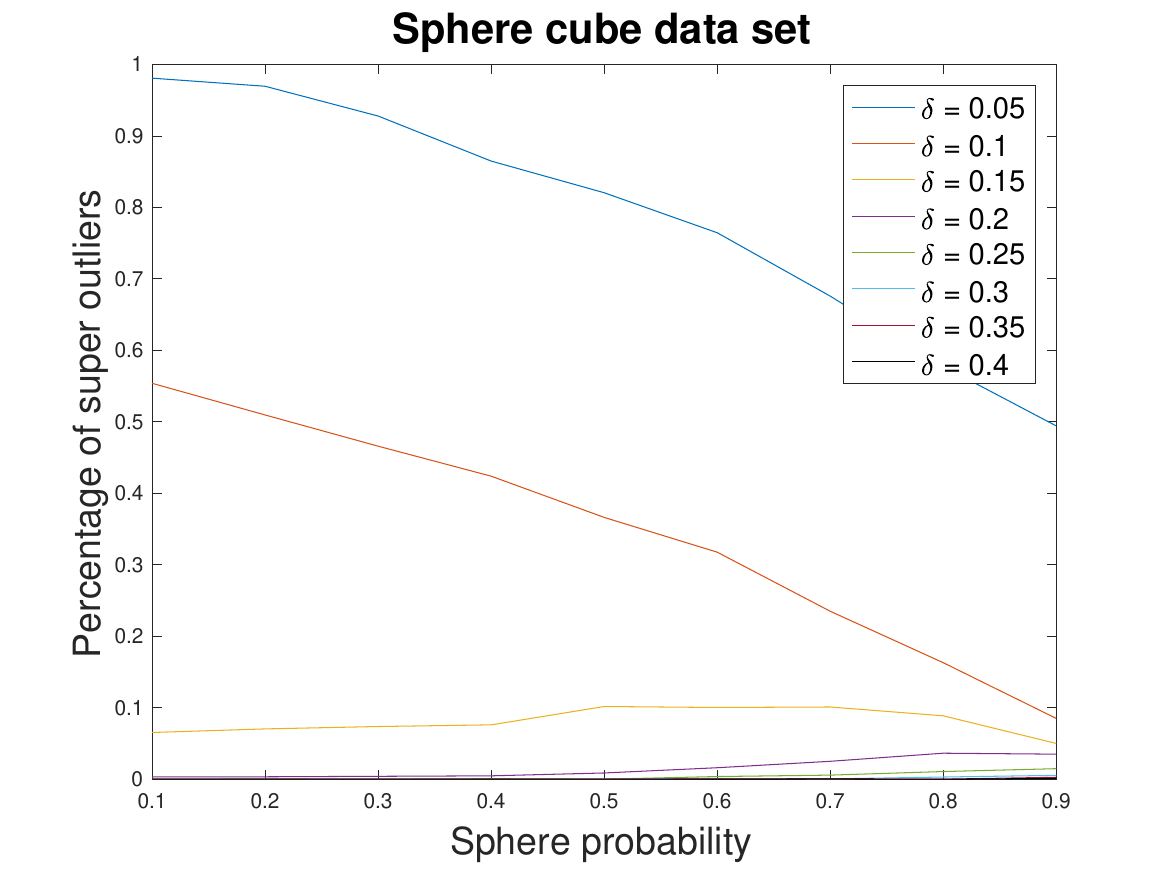}
\includegraphics[width=.49\textwidth]{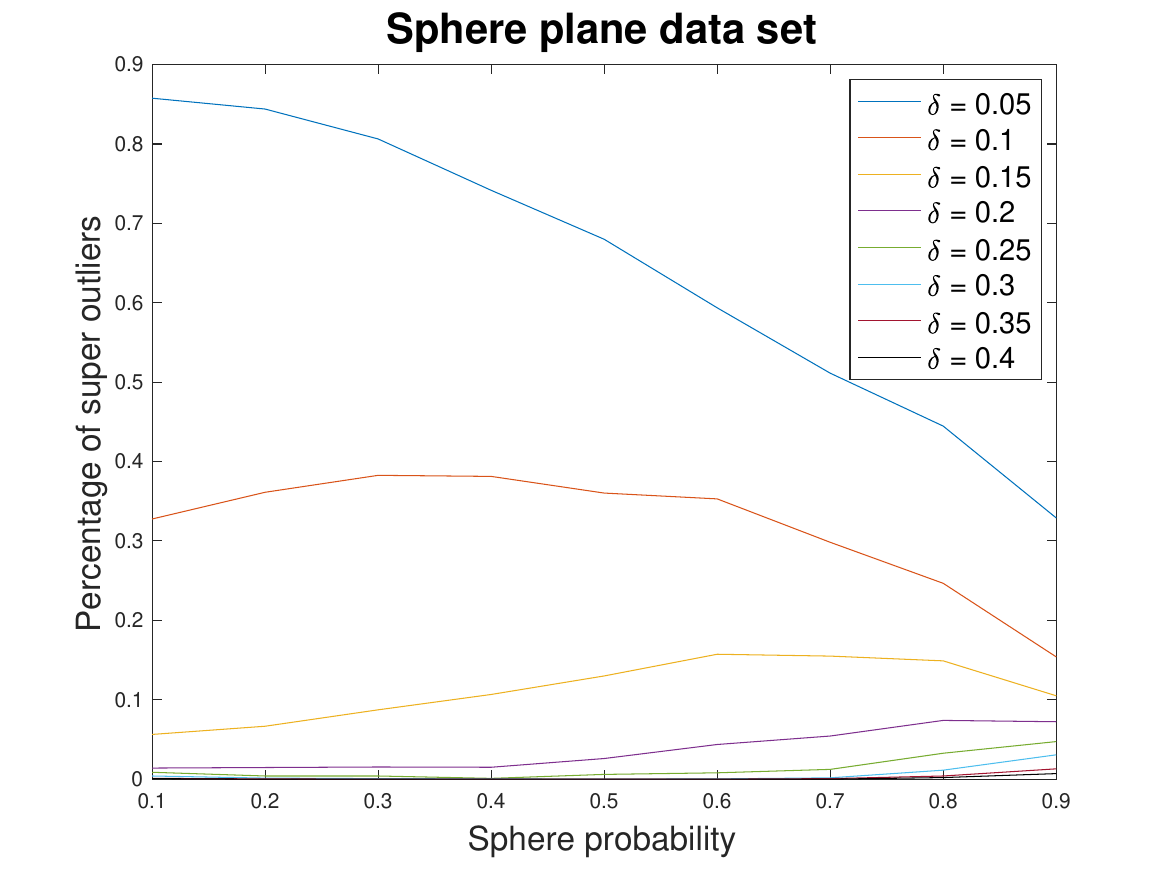}\\
\includegraphics[width=.49\textwidth]{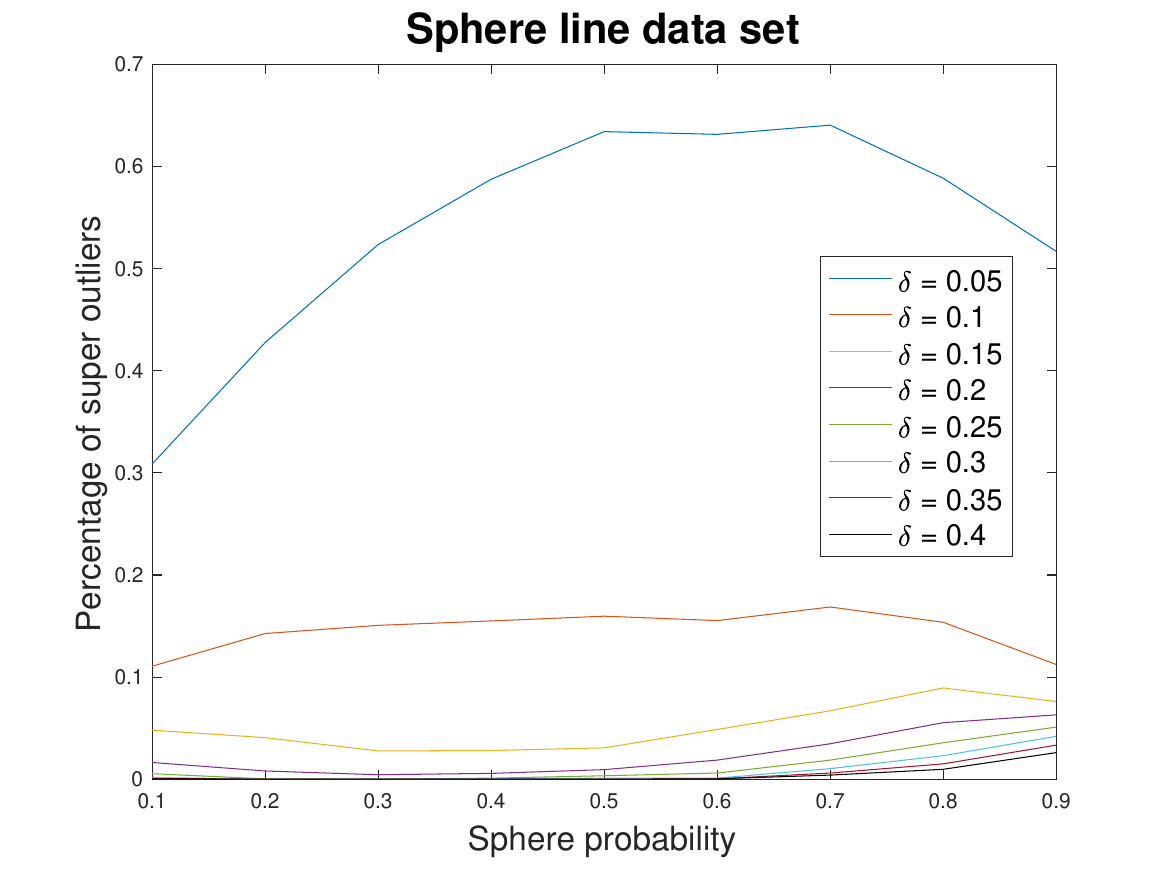}
\includegraphics[width=.49\textwidth]{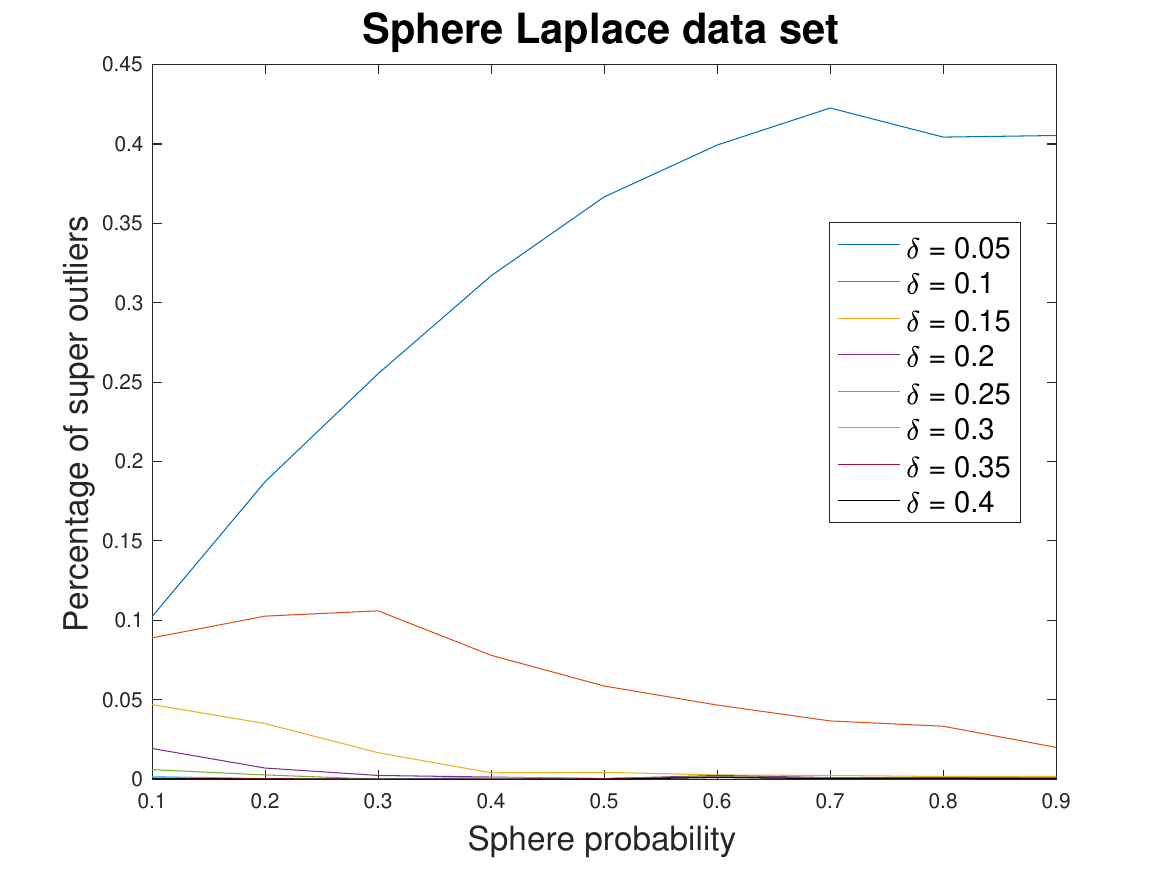}\\
\includegraphics[width=.49\textwidth]{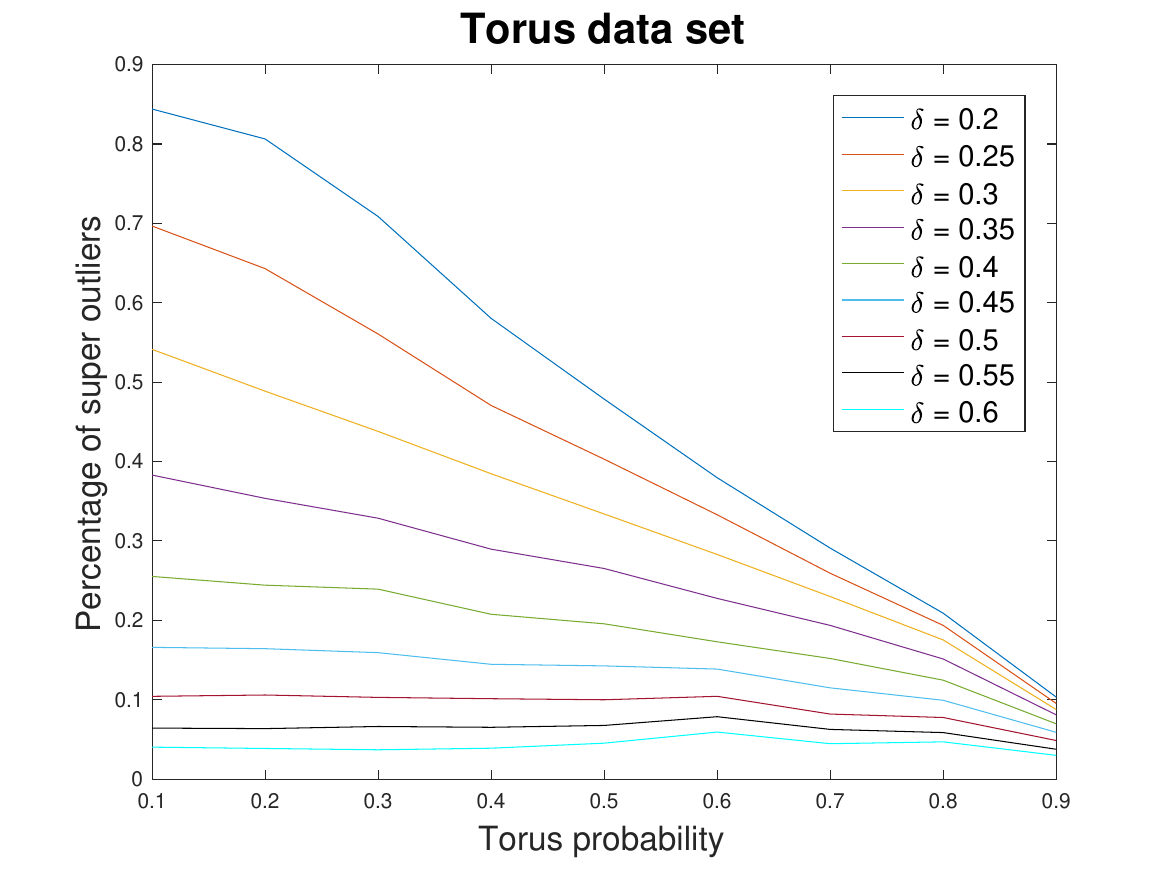}
\includegraphics[width=.49\textwidth]{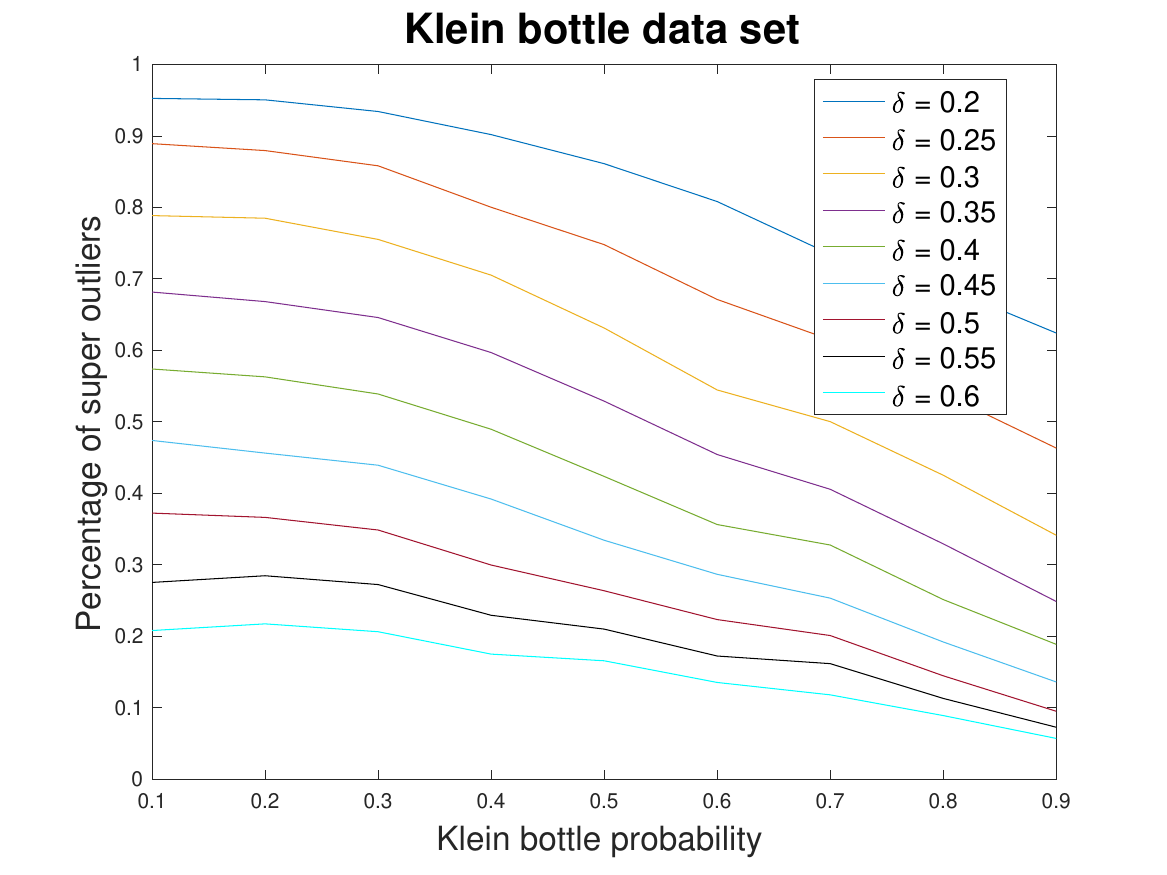}
\caption{\label{DeltaPlots} Influence of the choice of $\delta$ on the number of super outliers for different datasets.}
\end{figure}

Finally, we show how the number of super outliers depends on the choice of $\delta$ in Fig.~\ref{DeltaPlots}. We observe that small $\delta$ values lead to a drastic increase in the proportion of super outliers in the dataset. This underlines that, even though one can think of $\delta$ as a resolution parameter, depending on the dataset, the proportion of super outliers is an important factor to consider.

Overall, PH landmarks outperform all standard techniques on most datasets as well as two dense core subsets on the sphere-Laplace for a broad range of $\delta$-values. To obtain signal fractions that are high and stable over a long range of low sampling densities we recomment to choose $\delta$ as small as possible without having too many super outliers in the data. 

\clearpage

\section{Discussion}\label{subsec:LandmarkDiscussion}

We proposed novel outlier-robust landmark selection techniques: $k --$ landmarks and two versions of PH landmarks. PH landmark selection is the first landmark selection method developed specifically for the use of PH.

The $k --$ landmarks outperformed existing standard landmark selection methods -- random selection and  maxmin selection -- in many cases. They however tended to do so for large sampling densities, which are not very relevant for data subsampling.
While $k --$ landmarks did meet our goals for an outlier-robust landmark selection technique, this algorithm has a high computational cost and is difficult to use on a dataset with unknown properties as one has to predetermine the number of outliers for the algorithm to find.

We found that both versions of PH landmarks outperformed the existing standard landmark selection techniques on datasets containing noise, in particular for low landmark sampling densities. We observed this for a wide range of $\delta$ resolution values.
In most of our datasets, PH landmarks II (vital landmarks) restricted to dimension 1 for calculating the outlierness values slightly outperformed the PH landmarks I (vital landmarks) considering all dimensions (which coincides with a restriction to dimension 0 in our cases) for datasets with a high signal content and a low noise content, while PH landmarks I performed better for higher noise content. 

Dimension 1 PH landmarks II further outperformed dense core subsets showing that the method can capture richer information than just considering the $K$-th nearest neighbour as shown on the sphere-Laplace dataset. On our other datasets, both versions of PH landmarks were outperformed by dense core subsets due to the fact that in these datasets the signal points tend to be located in denser parts of the data than noise points. 
Overall, however, we consider PH landmarks to be a more practical approach than dense core subsets since they only require the choice of one parameter, $\delta$. In contrast, when applying dense core subsets as suggested in~\citep{Adams2014}, one needs to choose the parameter $K$ to obtain a density estimate, followed by the number of densest points to be chosen from the dataset, which are then used to obtain  maxmin landmarks. In~\citep{deSilva2004}, the authors observed that when studying dense core subsets the choices of parameters can indeed result in strong topological differences in the selected point clouds. With only one parameter choice we expect both versions of PH landmarks to deliver more consistent results.

In practice, we hypothesise that PH landmarks II using $out_{\text{PH}}^{1}(y)$ will lead to better results in high-dimensional data, especially in cases where signal data is sparse. The parameter $\delta$ would then have to be large enough to ensure that local PH in dimension 1 produces non-trivial barcodes for a large proportion of data points. For very large datasets with potentially high noise content, calculating $out_{\text{PH}}^{0,1,2}(y)$ or even a version restricted to dimension 0 for PH landmarks I could therefore be computationally advantageous.

We believe that PH landmarks contribute a valuable alternative to the current standard landmark selection techniques for PH, in particular for noisy datasets. The PH outlierness values as well as the number of super outliers with respect to the $\delta$ resolution could further be used to provide interesting insight into datasets for exploratory data analysis and could be incorporated in novel data analysis techniques which aim to take properties of shape of the data into account without performing computationally expensive PH analyses of the full dataset. Identifying representative and vital landmarks in datasets could be of interest to study processes underlying the data, for example in biology. While methods from Machine Learning are excellent at classifying such data, the underlying mechanisms often remain hidden to these methods and results by themselves can be difficult to interpret.
Building on our work, it would be interesting to investigate different definitions of PH outlierness in future studies. Testing the applicability of PH landmarks to real-world datasets with different types of noise or different levels of sparsity will be the subject of future work.


\acks{The author is very grateful for feedback and input from Vidit Nanda, Jared Tanner, and Heather Harrington as well as helpful suggestions from Scott Balchin. BJS is a member of the Centre for Topological Data Analysis, funded by the EPSRC grant (EP/R018472/1), and acknowledges funding from the EPSRC and MRC (EP/G037280/1) as well as F. Hoffmann-La Roche AG during her doctoral studies. 
}


\newpage

\appendix
\section{Appendix}
\label{app:theorem}



\subsection{Pseudocode for  maxmin algorithm~\citep{Adams2014}}\label{sec:Pseudocodemaxmin}

\begin{algorithm}
\caption[The  maxmin algorithm~\citep{Adams2014}]{The  maxmin algorithm~\citep{Adams2014}}
\label{alg:maxmin}
\begin{algorithmic} 
\REQUIRE Data points $D = \{y_1,\dots,y_N\}$, \\ a distance function $d:D\times D \rightarrow \R$ \\ number of landmarks $m$.
    \ENSURE A set of $m$ maxmin landmarks $L = \{l_1,\dots,l_m\}$.
    \STATE Select $y \in D$ at random
    \STATE $l_1 \leftarrow y$
    \STATE $L _1\leftarrow \{ l_1\}$
     \STATE $D_1 \leftarrow D\setminus \{ l_1\}$
\FOR{i = 2 \TO m}
	\FORALL{$y \in D_{i-1}$}
		\STATE Calculate $d(y,L_{i-1})$
	\ENDFOR
	\STATE Find $l_i$ such that $d(l_i,L_{i-1}) = \max_{y\in D}{d(y,L_{i-1})}$
	\STATE $L_i \leftarrow L_{i-1} \cup \{ l_i\}$
	\STATE $D_i \leftarrow D_{i-1} \setminus \{ l_i\}$
\ENDFOR
\end{algorithmic}
\end{algorithm}

\newpage

\subsection{Pseudocode for $k$--means$--$ algorithm~\citep{Chawla2013}}\label{sec:Pseudocodekminusminus}

\begin{algorithm}
\caption{The $k$--means$--$ algorithm~\citep{Chawla2013}}
\label{alg:kminusminus}
\begin{algorithmic} 
    \REQUIRE Data points $D = \{y_1,\dots,y_N\}$, \\ a distance function $d:D\times D \rightarrow \R$, \\ number of clusters $k$ and number of outliers $j$.
    \ENSURE A set of $k$ cluster centers $\hat{L} = \{\hat{l}_1,\dots,\hat{l}_k\}$, \\ a set of $j$ outliers $O = \{o_1,\dots,o_j \}$, $O \subset D$.
    \STATE $\hat{L}_0 \leftarrow \{k \text{ random points of } D\}$
    \STATE $i \leftarrow 1$
    \WHILE{(No convergence achieved)}
            \FORALL{$y \in D$}  \STATE compute $d(y,\hat{L}_{i-1})$ \ENDFOR
            \STATE Re-order the points in $D$ such that $d(y_1, \hat{L}_{i-1}) \geq d(y_2, \hat{L}_{i-1}) \geq \dots \geq d(y_N, \hat{L}_{i-1}) $
            \STATE $O_i \leftarrow \{y_1,\dots,y_k\}$
            \STATE $D_i \leftarrow D\setminus O_i = \{y_{k+1},\dots,y_N\}$
            \FOR{r = 1 \TO k}
            	\STATE {$P_r \leftarrow \{y \in D_i \ | \ c(y, \hat{L}_{i-1}) = \hat{l}_{i-1,r} \}$}
		\STATE {$\hat{l}_{i,r} \leftarrow \text{mean}(P_r)$}
            \ENDFOR
           \STATE $\hat{L}_i \leftarrow \{\hat{l}_{i,1}, \dots,  \hat{l}_{i,k}\}$
           \STATE $i \leftarrow i + 1$
    \ENDWHILE
\end{algorithmic}
\end{algorithm}

\newpage

\vskip 0.2in
\bibliography{Landmarks}

\end{document}